\documentclass[11pt]{amsart}

\usepackage[T1]{fontenc}
\usepackage{lmodern}
\usepackage{amsmath,amssymb,amsthm}
\usepackage[letterpaper,textwidth=6.28in,textheight=8.8in,centering]{geometry}
\usepackage[hidelinks,hypertexnames=false]{hyperref}

\hypersetup{
  pdftitle={The dual Orlicz-Minkowski problem for capillary hypersurfaces},
  pdfauthor={Ya Gao},
  pdfsubject={The dual Orlicz-Minkowski problem for capillary hypersurfaces}
}

\newcommand{\Ct}{\mathcal{C}_\theta}
\newcommand{\Hh}{\overline{\mathbb R^{n+1}_+}}
\newcommand{\Sp}{\mathbb S^n_+}
\newcommand{\dd}{\mathrm d}
\newcommand{\tr}{\operatorname{tr}}
\newcommand{\gradR}{\nabla_{\mathbb R^{n+1}}}
\newcommand{\hessR}{\nabla^2_{\mathbb R^{n+1}}}

\theoremstyle{plain}
\newtheorem{theorem}{Theorem}[section]
\newtheorem{proposition}[theorem]{Proposition}
\newtheorem{lemma}[theorem]{Lemma}
\newtheorem{corollary}[theorem]{Corollary}
\theoremstyle{remark}
\newtheorem*{remark}{Remark}
\numberwithin{equation}{section}

\allowdisplaybreaks
\title[The dual Orlicz-Minkowski problem for capillary hypersurfaces]
{The dual Orlicz--Minkowski problem for capillary hypersurfaces}
\author{Ya Gao$^{\dagger,\ast}$}
\address{$^{\dagger}$School of Mathematics and Statistics,
South-Central Minzu University, Wuhan 430074, China.}
\email{Echo-gaoya@outlook.com}
\thanks{$\ast$ Corresponding author}
\date{}

\begin{document}

\begin{abstract}
We study curvature flows for the dual Orlicz--Minkowski problem for
strictly convex capillary hypersurfaces in the Euclidean half-space.
The prescribed function has the form $\phi(h)G(X)$, where
$X=\nabla h+h\nu$ and $G$ may be nonradial and nonhomogeneous. For acute contact angles, the evenness assumption together with
either of the two integral conditions yields a smooth
weighted-volume-preserving flow for all $t>0$ and a smooth solution
to the dual Orlicz--Minkowski equation.
Without the evenness assumption, two strict spherical-cap barriers
yield the same conclusions for an unnormalized flow, whose
limiting hypersurfaces solve the dual Orlicz--Minkowski equation
with $c=1$.
A uniform strict monotonicity condition under rescaling ensures
uniqueness for $c=1$ and convergence of the unnormalized flow
in any $C^k$ norm at an exponential rate. The problem considered here includes the capillary $L_p$-Minkowski
and $L_p$ dual Minkowski problems as special cases.
\end{abstract}

\maketitle

\section{Introduction}\label{sec:intro}

We consider the dual Orlicz-Minkowski problem for capillary hypersurfaces with a general position weight. Throughout this paper, let $n\ge2$, $\theta\in(0,\pi/2)$ and $\{E_i\}_{i=1}^{n+1}$ be the standard orthonormal basis of $\mathbb{R}^{n+1}$, $e:=-E_{n+1}$. Write
\[
 \mathbb R^{n+1}_+=\{y:\langle y,E_{n+1}\rangle>0\},\qquad
 \Ct=\{\xi\in\Hh:|\xi-\cos\theta\,e|=1\}.
\]
For a strictly convex capillary hypersurface, the inverse capillary
Gauss parametrization gives
$\nu(\xi)=\xi-\cos\theta\,e$ and $X=\nabla h+h\nu$, where $h$ is its
support function. Given positive smooth functions $f\in C^{\infty}(\Ct)$,
$\phi\in C^{\infty}((0,\infty))$ and $G\in C^{\infty}(\Hh\setminus\{0\})$, we consider the following Monge-Amp\`ere equation with Robin boundary condition
\begin{equation}\label{eq:problem}
 \begin{cases}
 c\phi(h)G(\nabla h+h\nu)\det(\nabla^2h+h\sigma)=f
       &\text{in }\Ct,\\
 \nabla_\mu h=\cot\theta\,h&\text{on }\partial\Ct,
 \end{cases}
\end{equation}
where the
metric $\sigma$ is the unit spherical metric, $\mu$ is the outward
conormal of $\partial\Ct$, and determinants are taken with respect
to $\sigma$. Here $h>0$, $\nabla^2h+h\sigma>0$, and $c>0$ is a constant. The capillary geometry and the measure formulation of
\eqref{eq:problem} are recalled in Section~\ref{sec:prelim}.

For the classical Minkowski problem and its regularity theory,
see~\cite{Ale38,Ale39,CY76,Nir53}.
Lutwak introduced the $L_p$ Brunn--Minkowski theory
in~\cite{Lut.JDG.38-1993.131,Lut.Adv.118-1996.244}; further developments
of the $L_p$, logarithmic and centroaffine Minkowski problems include
\cite{BLYZ.JAMS.26-2013.831,CW.Adv.205-2006.33,HLYZ.DCG.33-2005.699,LYZ.TAMS.356-2004.4359}.
Dual curvature measures and their associated Minkowski problems
were introduced by Huang--Lutwak--Yang--Zhang~\cite{HLYZ.Acta.216-2016.325}.
Further results on the dual Minkowski problem can be found
in~\cite{BHP.JDG.109-2018.411,HP.Adv.323-2018.114,Zha.CVPDE.56-2017.18,Zha.JDG.110-2018.543}.
For closed hypersurfaces, related $L_p$ dual existence, variational
and curvature-flow results can be found
in~\cite{BF.JDE.266-2019.7980,CHZ.MA.373-2019.953,CCL,CLFlow,HZ.Adv.332-2018.57,LSW.JEMSJ.22-2020.893}.

Equation~\eqref{eq:problem} is the capillary counterpart of the dual
Orlicz--Minkowski equation. In the closed setting, Chen-Liu-Lu-Xiang~\cite{CLLX} obtained smooth even solutions by a Gauss curvature
flow preserving a weighted volume; Liu-Lu~\cite{LiuLu} treated
this problem without symmetry under scale conditions. Taking
$\phi(s)=s^{1-p}$ and $G(y)=|y|^{q-n-1}$ gives the $L_p$ dual
Minkowski equation, which includes the $L_p$-Minkowski equation when
$q=n+1$; see~\cite{LYZ}.
Further curvature-flow approaches to Minkowski-type problems
are given in~\cite{BIS2019,BISOrlicz,CTWX2022,LSYY2022}.

For capillary hypersurfaces, the Minkowski and $L_p$-Minkowski problems
were studied by Mei-Wang-Weng~\cite{MWWMinkowski,MWWLp}.
Related capillary Minkowski-type problems have also been studied
in~\cite{HICurvature,HISChristoffel,MWWChristoffel,WZOrlicz}.
Curvature-flow approaches and related rigidity results were developed
in~\cite{HHI,LiLiu,MWWFlow,MWWSolitons}.

 In particular, Hu-Hu-Ivaki~\cite{HHI}
studied capillary $L_p$ Minkowski flows, and Li-Liu~\cite{LiLiu}
considered an Orlicz coefficient depending on the normal direction
and $h/l$, where
\[
 l(\xi)=\sin^2\theta+\cos\theta\langle\xi,e\rangle
\]
is the support function of the unit cap. The capillary $L_p$ dual
Minkowski problem for $p>q$ is treated in~\cite{GaoEarlier,Gao,HHY}.
This paper concerns general $\phi(h)G(X)$ and the associated
flow asymptotics. Unlike a prescribed function w.r.t. $(\xi,h/l)$,
the position weight $G(X)$ also depends on $\nabla h$; it introduces
additional terms in both the interior and boundary curvature estimates.

\subsection{The curvature flow and the main results}
To solve~\eqref{eq:problem}, we balance a curvature-dependent inward
velocity with the rescaling term $+X$. Let $M$ be a compact smooth
$n$-dimensional manifold with boundary, and
$X_0:M\to\Hh$ be an embedding of a smooth strictly convex capillary hypersurface
$\Sigma_0$ with positive support function. Consider a family of smooth, strictly convex capillary hypersurfaces $\Sigma_t:= X(M, t)$ satisfying
\begin{equation}\label{eq:flow}
 \begin{cases}
 \displaystyle
 \partial_t X=-\eta(t)
 \frac{(f\circ\widetilde\nu)\langle X,\nu\rangle K}
 {(1+\cos\theta\langle\nu,e\rangle)
       \phi(\langle X,\nu\rangle)G(X)}\,\widetilde\nu+X
       &\text{in }M\times(0,T),\\[4pt]
 \langle\widetilde\nu,e\rangle=0
       &\text{on }\partial M\times(0,T),\\
 X(\cdot,0)=X_0&\text{in }M,
 \end{cases}
\end{equation}
where $\nu$ and $K$ are its outer unit normal
and Gauss curvature, $\widetilde\nu=\nu+\cos\theta\,e$ is the
capillary Gauss map, and $\eta(t)>0$ is to be chosen. At a
hypersurface the normal velocity vanishes precisely when
$\eta fK=\phi(h)G(X)$, giving~\eqref{eq:problem} with $c=\eta^{-1}$.

As usual,
\emph{evenness} means invariance under horizontal reflection
$(y',y_{n+1})\mapsto(-y',y_{n+1})$. We first obtain solutions of \eqref{eq:problem} without prescribing the constant $c$. To this end, we choose $\eta(t)$ to preserve a suitable $G$-weighted volume; see Section~\ref{subsec:volume} for its precise definition. We obtain the following result.

\begin{theorem}\label{thm:even}
Let $\theta\in(0,\pi/2)$, $\Sigma_0\subset\Hh$ be an even,
smooth, strictly convex capillary hypersurface with positive support
function. Let $\phi\in C^\infty((0,\infty))$ be positive, 
$f\in C^\infty(\Ct)$ and $G\in C^\infty(\Hh\setminus\{0\})$ be
positive and even. Suppose that $\phi$ and $G$ satisfy either
\begin{equation}\tag{A+}\label{ass:plus}
 \int_0^1\frac{\dd s}{\phi(s)}<\infty,\qquad
 \int_1^\infty\frac{\dd s}{\phi(s)}=\infty,\qquad
 \int_{B_1\cap\Hh}G(y)\,\dd y<\infty;
\end{equation}
or
\begin{equation}\tag{A-}\label{ass:minus}
 \int_1^\infty\frac{\dd s}{\phi(s)}<\infty,\qquad
 \int_{B_1\cap\Hh}G(y)\,\dd y=\infty,\qquad
 \int_{\Hh\setminus B_1}G(y)\,\dd y<\infty.
\end{equation}
Then flow~\eqref{eq:flow}, with $\eta$
given by~\eqref{eq:eta}, has an even, smooth, strictly convex capillary
solution for all $t>0$. A subsequence of $\Sigma_t$ converges in
$C^\infty$ to an even, smooth, strictly convex capillary hypersurface
whose positive support function solves~\eqref{eq:problem} for some
constant $c$. The corresponding $G$-weighted volume is preserved
along the flow and at the limit.
\end{theorem}

\begin{remark}
The constant $c$ in Theorem~\ref{thm:even} is not prescribed
in advance but is determined by the limiting hypersurface.
More precisely,
\[
c=
\frac{\displaystyle\int_{C_\theta}
\frac{fh_\infty}{\phi(h_\infty)}\,d\sigma}
{\displaystyle\int_{C_\theta}
\frac{h_\infty G(X_\infty)}{K_\infty}\,d\sigma},
\]
where $h_\infty$, $X_\infty$ and $K_\infty$ denote the
support function, position vector and Gauss curvature
of the limit hypersurface, respectively. This formula follows from
the choice of $\eta(t)$ in Section~\ref{subsec:volume}.
\end{remark}

To obtain solutions of \eqref{eq:problem} with $c=1$, we instead take $\eta(t)\equiv1$ in \eqref{eq:flow}. In this case, we use spherical caps as barriers to obtain uniform bounds without imposing the evenness assumption. This leads to the following long-time existence and convergence result.

\begin{theorem}\label{thm:barrier}
Let $\theta\in(0,\pi/2)$ and $\Sigma_0=s_0\Ct$, where $s_0>0$.
Let $f\in C^\infty(\Ct)$, $\phi\in C^\infty((0,\infty))$ and
$G\in C^\infty(\Hh\setminus\{0\})$ be positive. Suppose that there
exist $0<s_-<s_0<s_+$ such that
\begin{equation}\label{ass:barriers}
 s_-^n\phi(s_-l(\xi))G(s_-\xi)
       >f(\xi)>
 s_+^n\phi(s_+l(\xi))G(s_+\xi)
       \qquad\text{for any }\xi\in\Ct.
\end{equation}
Then flow~\eqref{eq:flow} with $\eta\equiv1$ has a unique smooth
strictly convex capillary solution for all $t>0$. Its support function
satisfies $s_-l<h(\cdot,t)<s_+l$. A subsequence of $\Sigma_t$ converges
in $C^\infty$ to a smooth strictly convex capillary hypersurface whose
positive support function solves~\eqref{eq:problem} with $c=1$.
\end{theorem}

\begin{remark}
Condition~\eqref{ass:barriers} ensures that
$s_-l$ and $s_+l$ are a strict subsolution
and a strict supersolution, respectively,
of the support-function flow with $\eta=1$.
By the comparison principle, the solution
remains between these two barriers for all
time; see Lemma~\ref{lem:barriers}.
Unlike the integral conditions in
Theorem~\ref{thm:even},
\eqref{ass:barriers} requires neither
evenness nor integrability of $G$
near the origin or at infinity.
\end{remark}

The following scale limits provide such barriers for any positive
smooth function $f$ on $\Ct$.
\begin{corollary}\label{cor:scales}
Let $\theta\in(0,\pi/2)$, $\phi\in C^\infty((0,\infty))$ and
$G\in C^\infty(\Hh\setminus\{0\})$ be positive. 
Suppose that
\begin{equation}\label{ass:limits}
 \lim_{s\to 0^+}\min_{\xi\in\Ct}
       s^n\phi(sl(\xi))G(s\xi)=\infty,\qquad
 \lim_{s\to\infty}\max_{\xi\in\Ct}
       s^n\phi(sl(\xi))G(s\xi)=0.
\end{equation}
Then, for any positive $f\in C^\infty(\Ct)$,
the equation \eqref{eq:problem} admits a positive smooth strictly convex solution
with $c=1$.
\end{corollary}

We next establish the uniqueness of solutions to
\eqref{eq:problem} with $c=1$ and the convergence
of the flow under a condition that ensures
strict monotonicity of the left-hand side of
\eqref{eq:problem} under rescaling.

\begin{theorem}\label{thm:unique}
Let $n\ge2$, $\theta\in(0,\pi/2)$, 
$f\in C^\infty(C_\theta)$,
$\phi\in C^\infty((0,\infty))$ and
$G\in C^\infty(\Hh\setminus\{0\})$ be positive
functions. Suppose that there exists a positive constant
$\delta_*$ such that
\begin{equation}\tag{U}\label{ass:U}
n+\frac{s\phi'(s)}{\phi(s)}
+\left\langle
\nabla_{\mathbb R^{n+1}}\log G(X),X
\right\rangle
\le -\delta_*
\end{equation}
for all $s>0$ and
$X\in\overline{\mathbb R^{n+1}_+}\setminus\{0\}$.

Then, for any initial capillary spherical cap
$\Sigma_0=s_0C_\theta$ with $s_0>0$, the flow
\eqref{eq:flow} with $\eta(t)\equiv1$ admits a unique
smooth, strictly convex capillary solution for all
$t>0$. Moreover, the family converges smoothly
to a strictly convex capillary hypersurface whose
support function is the unique positive smooth strictly
convex solution of \eqref{eq:problem} with $c=1$.
The convergence is exponential in any $C^k$ norm,
as specified in \eqref{eq:rates}.
\end{theorem}

\begin{remark}
\leavevmode\par\nobreak
\noindent\textup{(a)}\enspace
Condition~\eqref{ass:U} ensures that
$\phi(h)G(X)\det A$ decreases at a uniform
logarithmic rate as the rescaling
factor increases.
Indeed, for any $s>0$ and
$X\in\Hh\setminus\{0\}$, it implies
\[
\frac{d}{d\log a}
\log\big(a^n\phi(as)G(aX)\big)
\le-\delta_*,
\qquad a>0.
\]
Thus $a^n\phi(as)G(aX)$ is strictly decreasing
in $a$, with a uniform logarithmic rate.

This condition has two consequences. First, it
implies the scale limits \eqref{ass:limits}, and
therefore guarantees the existence of spherical-cap
barriers for any positive smooth function $f$.
Second, the strict monotonicity under rescaling
allows the comparison principle to establish
uniqueness of the solution.
Furthermore, condition~\eqref{ass:U} also plays
an important role in the convergence of the flow.
Indeed, together with the a priori estimates,
it ensures that the evolution equation for
$w=\partial_t\log h$ has a uniformly negative
zeroth-order coefficient. The maximum principle
then implies the exponential decay of $w$.
Combined with the uniform higher-order estimates,
this yields smooth convergence of the flow;
see Section~\ref{sec:convergence} for details.
\par\smallskip\noindent
\textup{(b)}\enspace
Theorems~\ref{thm:even}, \ref{thm:barrier}
and~\ref{thm:unique} concern different choices of the
time-dependent factor $\eta(t)$ and different
assumptions on the prescribed function:
In Theorem~\ref{thm:even}, $\eta(t)$ is chosen to
preserve a suitable $G$-weighted volume. Under
the integral conditions \eqref{ass:plus} or
\eqref{ass:minus} and the evenness assumption,
the flow has a smooth subsequence
solving \eqref{eq:problem} for some constant
$c>0$, which is not prescribed in advance;
In Theorem~\ref{thm:barrier}, we take
$\eta(t)\equiv1$ and use spherical-cap barriers
instead of a preserved weighted volume. This
removes the evenness assumption and gives a subsequence solving
\eqref{eq:problem} with $c=1$;
Theorem~\ref{thm:unique} strengthens
Theorem~\ref{thm:barrier} under the additional
condition \eqref{ass:U}: the barriers exist for
every positive smooth $f$, the flow converges exponentially to the unique solution. 

\end{remark}

\subsection{Outline of the proof}

We first establish uniform bounds for the support
function. For the flow preserving a $G$-weighted
volume, we use the monotonicity formula and the
evenness assumption, following \cite{CLLX}.
For the flow with $\eta\equiv1$, the bounds follow
from comparison with spherical caps. The gradient
estimate then follows from the convexity of the
evolving hypersurfaces.

The main step is the curvature estimate. Following
\cite{HHI,LiLiu}, we first derive upper and lower
bounds for the Gauss curvature. We then adapt the
auxiliary function $Q$
from \cite[Section 4]{Gao} to bound the principal
radii of curvature. The main difficulty is to
control the additional terms arising from the
general position weight $G(X)$ in the interior and boundary estimates.
These estimates imply long-time existence of the
flows. The monotonicity formulas yield smooth subsequences, while condition
\eqref{ass:U} ensures uniqueness and exponential
convergence of the flow.

\section{Preliminaries}\label{sec:prelim}

\subsection{Capillary geometry}
We follow the notation and support-function description in
\cite[Section 2]{Gao}. A capillary hypersurface $\Sigma$ meets
$\partial\mathbb R^{n+1}_+$ at angle $\theta$, so its outer normal
satisfies $\langle\nu,e\rangle=-\cos\theta$ along $\partial\Sigma$.
Write $\widehat\Sigma$ for the convex body enclosed by $\Sigma$ and $\partial\mathbb{R}^{n+1}_{+}$. We consider bodies for which the origin lies in the
relative interior of the base. The inverse capillary Gauss map
$X=\widetilde\nu^{-1}:\Ct\to\Sigma$ gives
\begin{equation}\label{eq:geometry}
 \begin{gathered}
 h=\langle X,\nu\rangle>0,\qquad
 X=\nabla h+h\nu,\qquad A_{ij}=h_{ij}+h\sigma_{ij}>0,\\
 K=(\det A)^{-1},\qquad
 l=\sin^2\theta+\cos\theta\langle\xi,e\rangle,\qquad u=h/l.
 \end{gathered}
\end{equation}
The eigenvalues of $A$ with respect to $\sigma$ are the principal
radii. Conversely, a positive strictly convex support function
satisfying the Robin condition determines such a capillary body;
see~\cite[Section 2]{Gao} and~\cite[Proposition 2.4]{MWWMinkowski}.
All covariant derivatives below are taken in $\sigma$. We use local
orthonormal frames, write $(A^{ij})=(A_{ij})^{-1}$, and sum repeated
indices. At the boundary $e_n=\mu$ and Greek indices are tangential.

The identities needed below are
\begin{equation}\label{eq:l-u}
 \begin{gathered}
 \nabla^2l+l\sigma=\sigma,\qquad
 1-\cos\theta\le l\le\sin^2\theta,\\
 A=l\nabla^2u+u\sigma+\nabla l\otimes\nabla u+\nabla u\otimes\nabla l,
 \qquad X=l\nabla u+u\xi,
 \end{gathered}
\end{equation}
on $\partial\Ct$, we have
\begin{equation}\label{eq:boundary}
 \begin{gathered}
 h_\mu=\cot\theta\,h,\qquad l_\mu=\cot\theta\,l,\qquad u_\mu=0,\\
 h_{\alpha\mu}=A_{\alpha\mu}=0,\qquad
 u_{\alpha\mu}=-\cot\theta\,u_\alpha.
 \end{gathered}
\end{equation}
The mixed boundary identities are given in
\cite[Proposition 2.1]{Gao}.
We shall also use the geometric estimate from
\cite[Lemma 3.2]{Gao}:
\begin{equation}\label{eq:geometric-bounds}
 \max_{\Ct}|X|\le\frac{\max_{\Ct}h}{\sin\theta},\qquad
 |\nabla h|\le |X|,\qquad
 B_{\min h}\cap\Hh\subset\widehat\Sigma
       \subset B_{\max|X|}\cap\Hh.
\end{equation}
Here $B_s$ is the ball of radius $s$ centered at the origin. These
statements use strict convexity and the Robin condition, not the
particular prescribed equation.

In the smooth case, the corresponding Orlicz curvature measure can
be defined by
\[
 \widetilde C^c_{\phi,G}(\widehat\Sigma,\omega)
 =\frac1{n+1}\int_{X(\omega)}
 (1+\cos\theta\langle\nu,e\rangle)
       \phi(\langle X,\nu\rangle)G(X)\,\dd S_\Sigma,
 \qquad\omega\subset\Ct.
\]
The area formula $X^*(\dd S_\Sigma)=K^{-1}\dd\sigma$ gives
\begin{equation}\label{eq:measure}
 \dd\widetilde C^c_{\phi,G}
       =\frac{l}{n+1}\phi(h)G(X)K^{-1}\dd\sigma.
\end{equation}
Thus the first equation in~\eqref{eq:problem} is equivalent to
$c\,\dd\widetilde C^c_{\phi,G}=fl\,\dd\sigma/(n+1)$.
For power weights this agrees with the measure formulation in
\cite[Section 2]{Gao}.

\subsection{Weighted volumes and short-time existence}
\label{subsec:volume}

For the weighted-volume-preserving flow in
Theorem~\ref{thm:even}, we use the interior weighted volume
under~\eqref{ass:plus} and the exterior weighted volume
under~\eqref{ass:minus}:
\begin{equation}\label{eq:volumes}
 V_{G,+}(\widehat\Sigma)
   =\int_{\widehat\Sigma}G(y)\,\dd y,
 \qquad
 V_{G,-}(\widehat\Sigma)
   =\int_{\Hh\setminus\widehat\Sigma}G(y)\,\dd y.
\end{equation}
The corresponding weighted volume is finite under
the respective assumptions.

Let $\rho$ be the radial function of $\widehat\Sigma$ on the closed
upper hemisphere $\Sp$. At corresponding radial and Gauss points,
the standard radial identities give
\begin{equation}\label{eq:radial}
 X(\xi,t)=\rho(\omega,t)\omega,\qquad
 \dd\omega=\frac{h}{K\rho^{n+1}}\,\dd\sigma,\qquad
 \frac{\rho_t}{\rho}=\frac{h_t}{h}.
\end{equation}
The derivatives in the last identity are at fixed $\omega$ and fixed
$\xi$, respectively; see~\cite[Section 2]{Gao} for the radial maps
and~\cite[Section 2]{CLLX} for this variation formula.
In particular, for a fixed $a>0$,
\begin{equation}\label{eq:first-variation}
 \frac{\dd}{\dd t}\int_{\Sp}\int_a^{\rho(\omega,t)}
                G(s\omega)s^n\,\dd s\,\dd\omega
       =\int_{\Ct}\frac{G(X)}K h_t\,\dd\sigma.
\end{equation}
Formula~\eqref{eq:first-variation}
also gives the derivatives of $V_{G,+}$ and $V_{G,-}$, with signs $+$
and $-$, by separating a fixed inner ball or exterior domain. Thus it
remains valid when $G$ is singular at the origin.

Taking the normal component of~\eqref{eq:flow} in Gauss coordinates, we have
\begin{equation}\label{eq:scalar}
 \begin{cases}
 \displaystyle h_t=h-\frac{\eta(t)fhK}{\phi(h)G(X)}
       &\text{in }\Ct\times(0,T),\\
 h_\mu=\cot\theta\,h&\text{on }\partial\Ct\times(0,T),\\
 h(\cdot,0)=h_0&\text{in }\Ct.
 \end{cases}
\end{equation}
The scalar equation and the geometric flow are equivalent up to
boundary-preserving tangential reparametrization. To preserve the
appropriate volume in~\eqref{eq:volumes}, its variation must vanish.
Substituting the first term of~\eqref{eq:scalar} into~\eqref{eq:first-variation} therefore
uniquely determines
\begin{equation}\label{eq:eta}
 \eta[h]=\frac{N[h]}{D[h]},\qquad
 N[h]=\int_{\Ct}\frac{hG(X)}K\,\dd\sigma,\qquad
 D[h]=\int_{\Ct}\frac{fh}{\phi(h)}\,\dd\sigma.
\end{equation}
Both integrals are positive. 

For $\eta\equiv1$, any positive strictly convex stationary solution
of~\eqref{eq:scalar} solves~\eqref{eq:problem} with $c=1$.
Let $h_\varepsilon=h+\varepsilon v$, where
$v_\mu=\cot\theta\,v$ on $\partial\Ct$.
Using~\eqref{eq:radial} and
$\left.\partial_\varepsilon\rho_\varepsilon
\right|_{\varepsilon=0}=\rho v/h$, we obtain
\begin{equation}\label{eq:eta-variation}
\begin{aligned}
 N'[h]v
 =\left.\frac{\dd}{\dd\varepsilon}
 N[h_\varepsilon]\right|_{\varepsilon=0}
 &=\int_{\Sp}
 \big((n+1)G(X)+\gradR G(X)\cdot X\big)
 \frac{\rho^{n+1}v}{h}\,\dd\omega\\
 &=\int_{\Ct}
 \big((n+1)G(X)+\gradR G(X)\cdot X\big)
 \frac vK\,\dd\sigma,\\[4pt]
 D'[h]v
 =\left.\frac{\dd}{\dd\varepsilon}
 D[h_\varepsilon]\right|_{\varepsilon=0}
 &=\int_{\Ct}\frac f{\phi(h)}
 \left(1-\frac{h\phi'(h)}{\phi(h)}\right)v\,\dd\sigma,\\[4pt]
 \eta'[h]v
 =\left.\frac{\dd}{\dd\varepsilon}
 \eta[h_\varepsilon]\right|_{\varepsilon=0}
 &=\frac{D[h]N'[h]v-N[h]D'[h]v}{D[h]^2}\\
 &=\frac{N'[h]v-\eta[h]D'[h]v}{D[h]}.
\end{aligned}
\end{equation}

\begin{lemma}[Short-time existence]\label{lem:short-time}
Let $h_0\in C^\infty(\Ct)$ satisfy $h_0>0$,
$\nabla^2h_0+h_0\sigma>0$ and
$(h_0)_\mu=\cot\theta\,h_0$ on $\partial\Ct$.
For $\eta$ given by~\eqref{eq:eta} or $\eta\equiv1$,
there exist $\alpha\in(0,1)$ and $T>0$ such that~\eqref{eq:scalar}
has a unique positive strictly convex solution
\[
 h\in C^{2+\alpha,1+\alpha/2}(\Ct\times[0,T])
       \cap C^\infty(\Ct\times(0,T]).
\]
If $f$, $G$ and $h_0$ are even, then $h(\cdot,t)$ remains even.
\end{lemma}
\begin{proof}
Set $h=lu$, so $u_\mu=0$, and denote the resulting right-hand side
by $\mathcal F(u)$. Near $u_0=h_0/l$ the equation is uniformly
parabolic. By~\eqref{eq:eta-variation}, its linearization is an
ordinary parabolic Neumann operator $L$ plus an integral term
\[
\begin{gathered}
 (Bv)(\xi,t)=a(\xi,t)\int_{\Ct}b(\zeta,t)v(\zeta,t)\,\dd\sigma(\zeta),
 \\
 \|Bv\|_{C^{\alpha,\alpha/2}}
 \le CT^{1-\alpha/2}\|v\|_{C^{2+\alpha,1+\alpha/2}}
 \quad\text{if }v(\cdot,0)=0.
\end{gathered}
\]
Here the parabolic H\"older norms of $a,b$ are bounded independently
of small $T$; $B=0$ if $\eta=1$.
Linear Neumann Schauder theory~\cite[Chapter V]{Lieberman}
and a Neumann series give invertibility for small $T$ on variations
with $v(\cdot,0)=0$ and $v_\mu=0$.
For the inverse-function argument, take the approximate extension
\[
 \bar u_t-\Delta_\sigma\bar u
 =\mathcal F(u_0)-\Delta_\sigma u_0,\qquad
 \bar u_\mu=0,\qquad \bar u(\cdot,0)=u_0.
\]
Its residual $\bar u_t-\mathcal F(\bar u)$ vanishes at $t=0$;
using exponents $0<\alpha<\alpha_1<1$, its
$C^{\alpha,\alpha/2}$ norm tends to zero with $T$.
The inverse-function theorem yields existence and uniqueness.
Positivity and strict convexity persist by continuity,
positive-time smoothness follows by parabolic bootstrapping,
and evenness follows from uniqueness.
\end{proof}

\subsection{Rescaling of curvature flow}\label{subsec:rescaling}
The construction above follows the weighted-volume approach
of~\cite{CLLX}. For homogeneous weights it also has the rescaling
interpretation used in~\cite[(1.2)--(1.3)]{HHI}.
To distinguish this from taking $\eta=1$, start with the following flow
\[
 \partial_tY=-\frac{(f\circ\widetilde\nu_Y)h_YK_Y}
 {l_Y\phi(h_Y)G(Y)}\,\widetilde\nu_Y,
\]
where $h_Y=\langle Y,\nu_Y\rangle$ and
$l_Y=1+\cos\theta\langle\nu_Y,e\rangle$, with the same capillary
boundary condition. Define
$X(\cdot,\tau(t)):=\lambda(t)Y(\cdot,t)$ and
$\tau(t)=\log\lambda(t)$, where $\lambda(0)=1$ and
$\dot\lambda:=\frac{d\lambda(t)}{dt}>0$.
Since $h=\lambda h_Y$ and $K=\lambda^{-n}K_Y$, we have
\begin{equation}\label{eq:rescaling}
 \partial_\tau X=X-
 \frac{\lambda^n}{\dot\tau}
 \frac{(f\circ\widetilde\nu)hK}
 {l\phi(h/\lambda)G(X/\lambda)}\,\widetilde\nu,
\end{equation}
where $\dot\tau:=\frac{d\tau(t)}{dt}=\frac{\dot\lambda(t)}{\lambda(t)}$. For general weights, $h/\lambda$ and $X/\lambda$
cannot be replaced by $h$ and $X$ using a scalar depending only on
time. Thus the two choices in~\eqref{eq:scalar} are not claimed to be
related by such a rescaling with fixed $\phi$ and $G$.

For comparison, suppose $\phi(s)=s^{1-p}$ and
$G(aY)=a^{q-n-1}G(Y)$ with $q\ne0$. Let $V_G$ denote the finite
interior volume for $q>0$ or exterior volume for $q<0$. Then
$V_G(a\widehat\Sigma)=a^qV_G(\widehat\Sigma)$. Along the
flow, the choice
\[
 \lambda(t)=\left(\frac{V_G(\widehat\Sigma_0)}
                         {V_G(\widehat M_t)}\right)^{1/q},
 \qquad
 \dot\tau=-\frac1q\frac{\dd}{\dd t}\log V_G(\widehat M_t)
          =\frac{D[h_Y]}{N[h_Y]}
\]
uses the body $\widehat M_t$ enclosed by $Y(M,t)$.
The scaling of $N$ and $D$ gives
$\dot\tau=\lambda^{q-p}/\eta[h]$, and~\eqref{eq:rescaling}
becomes~\eqref{eq:flow}, with time variable $\tau$ and
$\eta=N/D$. In particular, $G=1$ and $q=n+1$ give exactly the
volume rescaling in~\cite{HHI}. The coefficient choice in
\eqref{eq:eta}, rather than this homogeneity calculation, is used
for the general weights in Theorem~\ref{thm:even}.

\subsection{Differential identities}\label{subsec:dilation}
Define the logarithmic rescaling derivative
\begin{equation}\label{eq:E}
 E(s,Y)=n+\frac{s\phi'(s)}{\phi(s)}
               +\gradR\log G(Y)\cdot Y.
\end{equation}
With the above notation, condition~\eqref{ass:U}
is equivalent to
\[
E(s,Y)\le-\delta_*
\qquad
(s>0,\ Y\in\Hh\setminus\{0\}).
\]
Its significance follows from
\[
 \frac{\dd}{\dd\log a}\log\big(a^n\phi(as)G(aY)\big)=E(as,aY).
\]
Consequently,
\begin{equation}\label{eq:dilation-inequality}
 \begin{cases}
 a^n\phi(as)G(aY)\le a^{-\delta_*}\phi(s)G(Y),&a\ge1,\\
 a^n\phi(as)G(aY)\ge a^{-\delta_*}\phi(s)G(Y),&0<a\le1.
 \end{cases}
\end{equation}
Taking $(s,Y)=(l(\xi),\xi)$ proves~\eqref{ass:limits}, since $\Ct$
is compact and avoids the origin. The curvature
estimates below use only boundedness of $E$ and do not require its negativity.

For later use, put $r=fK/(\phi(h)G(X))$ and $v=\log h$.
For a variation $v+\varepsilon z$,
\begin{equation}\label{eq:log-linearization}
 \left.\frac{\dd}{\dd\varepsilon}\right|_{\varepsilon=0}
        \log r[e^{v+\varepsilon z}]
 =-hA^{ij}z_{ij}
  -(2A^{ij}h_i+h\gamma_j)z_j-E(h,X)z,
\end{equation}
where $\gamma:=\gradR\log G(X)$ and
$\gamma_j:=\langle\gamma,e_j\rangle$. 
Let $\delta$ denote the first variation with
respect to $\varepsilon$ at $\varepsilon=0$.
Indeed, $\delta h=hz$, $\delta X=zX+h\nabla z$ and
$\delta\log\det A=nz+hA^{ij}z_{ij}+2A^{ij}h_i z_j$.
This single linearization will be used for both the lower Gauss
curvature bound and the decay of the logarithmic speed.

The Codazzi and commutation identities, in the notation
$\sigma_1=\tr A$ and $S=\tr(A^{-1})$, are
\begin{equation}\label{eq:codazzi}
 A_{ij,k}=A_{ik,j},\qquad
 \Delta A_{ij}=(\sigma_1)_{ij}+nA_{ij}-\sigma_1\delta_{ij},
\end{equation}
see~\cite[Lemma 2.2]{Gao}. For a smooth function $g$,
\begin{equation}\label{eq:composition}
 \begin{aligned}
 D_iX&=A_{ij}e_j,\qquad
 D^2_{ij}X=A_{ij,k}e_k-A_{ij}\nu,\\
 \nabla_i[g(X)]&=A_{ij}\langle\gradR g,e_j\rangle,\\
 \Delta[g(X)]&=\sum_i\hessR g(Ae_i,Ae_i)
     +\langle(\gradR g)^\top,\nabla\sigma_1\rangle
     -\langle\gradR g,\nu\rangle\sigma_1,
 \end{aligned}
\end{equation}
where $D$ denotes the covariant derivative in $\mathbb{R}^{n+1}$, and
$\top$ denotes the tangential projection onto $\nu^\perp$.

\section{Monotonicity formulas and lower-order estimates}\label{sec:lower}

Let $h$ be the short-time solution of~\eqref{eq:scalar} on its maximal
interval $[0,T)$. We establish bounds for $h$, $\nabla h$ and $\eta$.
The two choices of $\eta$ require different arguments at this stage.

\subsection{The weighted-volume-preserving flow}
Assume the hypotheses of Theorem~\ref{thm:even}, and choose
$\eta=N/D$. In the respective cases~\eqref{ass:plus}
and~\eqref{ass:minus}, define
\[
 \Phi_+(s)=\int_0^s\frac{\dd a}{\phi(a)},\qquad
 \Phi_-(s)=\int_s^\infty\frac{\dd a}{\phi(a)},\qquad
 J_\pm(t)=\int_{\Ct}f\Phi_\pm(h)\,\dd\sigma.
\]
The function $\Phi_+$ is increasing and tends to infinity at infinity,
whereas $\Phi_-$ is decreasing and tends to zero there.

\begin{lemma}\label{lem:entropy}
The weighted volume in~\eqref{eq:volumes} is constant,
and we have
\begin{equation}\label{eq:entropy}
 \frac{\dd J_+}{\dd t}
       =-\frac1\eta\int_{\Ct}\frac{G(X)}{hK}h_t^2\,\dd\sigma\le0,
 \qquad
 \frac{\dd J_-}{\dd t}
       =\frac1\eta\int_{\Ct}\frac{G(X)}{hK}h_t^2\,\dd\sigma\ge0
\end{equation}
under~\eqref{ass:plus} and~\eqref{ass:minus}, respectively.
\end{lemma}
\begin{proof}
By~\eqref{eq:first-variation}, \eqref{eq:scalar} and~\eqref{eq:eta}, we have
\[
\begin{aligned}
 \frac{\dd}{\dd t}V_{G,\pm}(\widehat\Sigma_t)
 &=\pm\int_{\Ct}\frac{G(X)}K
       \left(h-\frac{\eta fhK}{\phi(h)G(X)}\right)\dd\sigma\\
 &=\pm\left(N[h]-\eta[h]D[h]\right)=0.
\end{aligned}
\]
The same scalar equation gives
\[
 \frac f{\phi(h)}
 =\frac{G(X)}{\eta hK}(h-h_t),\qquad
 \int_{\Ct}\frac{G(X)}K h_t\,\dd\sigma=0.
\]
Then we have
\[
\begin{aligned}
 J_+'(t)
 &=\int_{\Ct}\frac f{\phi(h)}h_t\,\dd\sigma\\
 &=\frac1\eta\int_{\Ct}\frac{G(X)}K h_t\,\dd\sigma
   -\frac1\eta\int_{\Ct}\frac{G(X)}{hK}h_t^2\,\dd\sigma\\
 &=-\frac1\eta\int_{\Ct}\frac{G(X)}{hK}h_t^2\,\dd\sigma\le0,
\end{aligned}
\]
similarly, 
\[
\begin{aligned}
 J_-'(t)
 &=-\int_{\Ct}\frac f{\phi(h)}h_t\,\dd\sigma
   =\frac1\eta\int_{\Ct}\frac{G(X)}{hK}h_t^2\,\dd\sigma\ge0.
\end{aligned}
\]
The calculation for $J_+$ applies under~\eqref{ass:plus},
and that for $J_-$ applies under~\eqref{ass:minus}.
The identities hold for $t>0$. After integration in time,
continuity at $t=0$ allows the initial time to be included.
\end{proof}

\begin{lemma}\label{lem:even-lower}
Under the hypotheses of Theorem~\ref{thm:even}, there is a constant
$C>1$, independent of $t$, such that
\begin{equation}\label{eq:low-bounds}
 C^{-1}\le h\le C,\qquad |\nabla h|\le C,\qquad
 C^{-1}\le\eta\le C.
\end{equation}
\end{lemma}
\begin{proof}
We follow the support-function argument of~\cite[Lemma 3]{CLLX},
using reflection across the supporting plane as in~\cite[Lemma 4.1]{HHI}.
Write
\[
 R(t)=\max_{y\in\widehat\Sigma_t}|y|,\qquad
 m(t)=\min_{\xi\in\Ct}h(\xi,t),\qquad
 f_{\min}=\min_{\Ct}f,\quad f_{\max}=\max_{\Ct}f,
\]
where $|y|$ is the Euclidean norm; $|E|$ denotes area for
a subset $E$ and $\operatorname{Vol}$ denotes Euclidean
volume in $\mathbb R^{n+1}$.

Let $\mathcal R(y',y_{n+1})=(y',-y_{n+1})$ and set
\[
 \Xi_t=\operatorname{conv}
       \bigl(\widehat\Sigma_t\cup\mathcal R(\widehat\Sigma_t)\bigr),
 \qquad h_{\Xi_t}(z)=\max_{y\in\Xi_t}\langle y,z\rangle,
 \quad z\in\mathbb S^n.
\]
By uniqueness in Lemma~\ref{lem:short-time}, horizontal evenness is
preserved, and hence $\Xi_t=-\Xi_t$. Since
$\nu_{n+1}(\xi)\ge\cos\theta>0$, reflection does not increase
the scalar product with $\nu(\xi)$, so
\[
 h_{\Xi_t}(\nu(\xi))=h(\xi,t),\qquad
 \max_{y\in\Xi_t}|y|=R(t).
\]
Choose $y_t\in\widehat\Sigma_t$ with $|y_t|=R(t)$ and let
$x_t=y_t/R(t)$. Then $\pm R(t)x_t\in\Xi_t$, and hence
\begin{equation}\label{eq:reflection-bounds}
 R(t)|\langle x_t,\nu(\xi)\rangle|
 \le h_{\Xi_t}(\nu(\xi))=h(\xi,t)\le R(t),
 \qquad \xi\in\Ct.
\end{equation}
If $h(\xi_t,t)=m(t)$, the same symmetry gives
\[
 \widehat\Sigma_t\subset\Xi_t\subset
 \{y:|\langle y,\nu(\xi_t)\rangle|\le m(t)\},\qquad
 \widehat\Sigma_t\subset B_{R(t)}.
\]

For $x\in\mathbb S^n$ and $0<\delta\le1/2$, define
\[
 \mathcal B_\delta(x)
 =\{\xi\in\Ct:|\langle x,\nu(\xi)\rangle|<\delta\}.
\]
The map $\xi\mapsto\nu(\xi)$ preserves spherical area. Rotational
invariance of the unit sphere therefore gives
\begin{equation}\label{eq:belt}
\begin{aligned}
 |\mathcal B_\delta(x)|
 &\le|\{z\in\mathbb S^n:|\langle x,z\rangle|<\delta\}|\\
 &=2|\mathbb S^{n-1}|\int_0^\delta(1-s^2)^{(n-2)/2}\,\dd s
 \le2|\mathbb S^{n-1}|\delta.
\end{aligned}
\end{equation}

\smallskip\noindent\emph{Case \eqref{ass:plus}.}
Choose $0<\delta_0\le1/2$ such that
$2|\mathbb S^{n-1}|\delta_0\le|\Ct|/2$.
By~\eqref{eq:entropy}, \eqref{eq:reflection-bounds} and~\eqref{eq:belt}, we have
\[
\begin{aligned}
 J_+(0)&\ge J_+(t)
 =\int_{\Ct}f\Phi_+(h)\,\dd\sigma\\
 &\ge f_{\min}
       \int_{\Ct\setminus\mathcal B_{\delta_0}(x_t)}
                      \Phi_+(h)\,\dd\sigma\\
 &\ge f_{\min}|\Ct\setminus\mathcal B_{\delta_0}(x_t)|
                      \Phi_+(\delta_0R(t))\\
 &\ge\tfrac12f_{\min}|\Ct|\Phi_+(\delta_0R(t)).
\end{aligned}
\]
Since $\Phi_+$ is strictly increasing and tends to infinity,
$R(t)\le R_0$ for a fixed $R_0>0$.
Let $V_0=V_{G,+}(\widehat\Sigma_0)>0$. By~\eqref{ass:plus},
$G\in L^1(B_{R_0}\cap\Hh)$, so there exists $v_0>0$ such that, for any measurable set $\Omega$,
\[
 \Omega\subset B_{R_0}\cap\Hh,\quad\operatorname{Vol}(\Omega)<v_0
 \quad\Longrightarrow\quad\int_\Omega G\,\dd y<V_0.
\]
Volume preservation implies $\operatorname{Vol}(\widehat\Sigma_t)\ge v_0$.
By the above inclusion, we have
\[
 v_0\le\operatorname{Vol}(\widehat\Sigma_t)
 \le2m(t)\,\omega_nR_0^n,
 \qquad
 m(t)\ge\frac{v_0}{2\omega_nR_0^n}>0,
\]
where $\omega_n$ is the volume of the unit ball in $\mathbb R^n$.

\smallskip\noindent\emph{Case \eqref{ass:minus}.}
Suppose $m(t_j)\to0$. Choose $\xi_j\in\Ct$ with
$h(\xi_j,t_j)=m(t_j)$ and, after taking a subsequence, let
$\nu(\xi_j)\to\nu_\infty\in\mathbb S^n$. For any
$y\in(B_1\cap\Hh)\setminus\nu_\infty^\perp$,
\[
 |\langle y,\nu(\xi_j)\rangle|
 \longrightarrow|\langle y,\nu_\infty\rangle|>0,
 \qquad m(t_j)\longrightarrow0.
\]
The above inclusion implies that
$y\notin\widehat\Sigma_{t_j}$ for all sufficiently large $j$. Writing $\boldsymbol1_E$ for the indicator of a set $E$,
we obtain $\boldsymbol1_{\Hh\setminus\widehat\Sigma_{t_j}}\to1$
almost everywhere on $B_1\cap\Hh$. Fatou's lemma and volume preservation yield
\[
\begin{aligned}
 \infty&=\int_{B_1\cap\Hh}G(y)\,\dd y\\
 &\le\liminf_{j\to\infty}\int_{B_1\cap\Hh}
     G(y)\boldsymbol1_{\Hh\setminus\widehat\Sigma_{t_j}}(y)\,\dd y\\
 &\le\liminf_{j\to\infty}V_{G,-}(\widehat\Sigma_{t_j})
 =V_{G,-}(\widehat\Sigma_0)<\infty,
\end{aligned}
\]
a contradiction. Hence $m(t)\ge m_0>0$.

Choose $0<\delta\le1/2$ such that
\[
 2f_{\max}|\mathbb S^{n-1}|\delta\Phi_-(m_0)
 \le\tfrac12J_-(0).
\]
The function $\Phi_-$ is decreasing. Splitting the integral into
$\mathcal B_\delta(x_t)$ and its complement gives
\[
\begin{aligned}
 J_-(0)&\le J_-(t)\\
 &=\int_{\mathcal B_\delta(x_t)}f\Phi_-(h)\,\dd\sigma
   +\int_{\Ct\setminus\mathcal B_\delta(x_t)}f\Phi_-(h)\,\dd\sigma\\
 &\le f_{\max}|\mathcal B_\delta(x_t)|\Phi_-(m_0)
   +f_{\max}|\Ct|\Phi_-(\delta R(t))\\
 &\le\tfrac12J_-(0)+f_{\max}|\Ct|\Phi_-(\delta R(t)).
\end{aligned}
\]
Thus $\Phi_-(\delta R(t))\ge J_-(0)/(2f_{\max}|\Ct|)>0$.
Since $\Phi_-(s)\to0$ as $s\to\infty$, this bounds $R(t)$ by a fixed $R_0$.

In both cases, \eqref{eq:geometric-bounds} gives fixed $m_0,R_0>0$ with
\[
 m_0\le h,\rho\le R_0,\qquad |\nabla h|\le|X|\le R_0.
\]
Finally, the radial identity~\eqref{eq:radial} gives
\[
 N[h]=\int_{\Sp}G(\rho\omega)\rho^{n+1}\,\dd\omega.
\]
Using $m_0\le h,\rho\le R_0$ and the definition of $D[h]$, we obtain
\[
\begin{aligned}
 |\Sp|\min_{[m_0,R_0]\times\Sp}
       \bigl(G(s\omega)s^{n+1}\bigr)
 &\le N[h]\le
 |\Sp|\max_{[m_0,R_0]\times\Sp}
       \bigl(G(s\omega)s^{n+1}\bigr),\\
 |\Ct|\min_{\Ct\times[m_0,R_0]}
       \frac{f(\xi)s}{\phi(s)}
 &\le D[h]\le
 |\Ct|\max_{\Ct\times[m_0,R_0]}
       \frac{f(\xi)s}{\phi(s)}.
\end{aligned}
\]
The four extrema are positive and finite by smoothness and compactness.
Hence there exist constants $0<c\le C$, independent of $t$, such that
$c\le N[h],D[h]\le C$. By~\eqref{eq:eta}, we have
\[
 \frac{c}{C}\le\eta[h]=\frac{N[h]}{D[h]}\le\frac{C}{c}.
\]
\end{proof}

\subsection{The flow with \texorpdfstring{$\eta=1$}{eta=1}}
Assume that~\eqref{ass:barriers} holds and $h_0=s_0l$, as in
Theorem~\ref{thm:barrier}. No symmetry assumption is made here.

\begin{lemma}\label{lem:barriers}
For all $t\in[0,T)$, the solution satisfies
\begin{equation}\label{eq:barrier-bounds}
 s_-<u=h/l<s_+.
\end{equation}
Then, there exists a constant $C>1$, independent of $t$,
such that
\[
 C^{-1}\le h,|X|,\rho\le C,\qquad |\nabla h|\le C.
\]
\end{lemma}
\begin{proof}
      Suppose, for a contradiction, that $u$ reaches $s_-$, and let
$t_0$ be the first time when this happens.
Since $u(\cdot,0)=s_0>s_-$, we have $t_0>0$.
There exists $\xi_0\in\Ct$ such that
\[
 u(\xi_0,t_0)=\min_{\Ct}u(\cdot,t_0)=s_-,
 \qquad u_t(\xi_0,t_0)\le0.
\]
If $\xi_0\in\partial\Ct$, the Neumann condition and~\eqref{eq:boundary} give
\[
 u_\alpha=u_\mu=0,\qquad u_{\alpha\mu}=-\cot\theta\,u_\alpha=0,
 \qquad (u_{\alpha\beta})\ge0,\quad u_{\mu\mu}\ge0.
\]
The same second-derivative inequalities hold at an interior minimum. Thus
\[
 u_t\le0,\qquad\nabla u=0,\qquad\nabla^2u\ge0.
\]
Hence~\eqref{eq:l-u} implies
\[
 h=s_-l,\qquad X=s_-\xi_0,\qquad
 A\ge s_-\sigma,\qquad K\le s_-^{-n}.
\]
Using~\eqref{eq:scalar} and~\eqref{ass:barriers}, we obtain
\[
 0\ge u_t
 =s_-\left(1-\frac{fK}{\phi(h)G(X)}\right)
 \ge s_-\left(1-
 \frac{f(\xi_0)}{s_-^n\phi(s_-l(\xi_0))G(s_-\xi_0)}\right)>0.
\]
At a first upper contact with $u=s_+$, the same identities give
\[
 u_t\ge0,\qquad\nabla u=0,\qquad\nabla^2u\le0,
 \qquad A\le s_+\sigma,\quad K\ge s_+^{-n},
\]
and therefore
\[
 0\le u_t
 \le s_+\left(1-
 \frac{f(\xi_0)}{s_+^n\phi(s_+l(\xi_0))G(s_+\xi_0)}\right)<0.
\]
Both contacts are impossible, so we have the inequality~\eqref{eq:barrier-bounds}.
Finally, by \eqref{eq:l-u} and~\eqref{eq:geometric-bounds}, we have
\[
 0<s_-\min_{\Ct}l<h<s_+\max_{\Ct}l,\qquad
 \min h\le\rho=|X|\le\frac{\max h}{\sin\theta},\qquad
 |\nabla h|\le|X|.
\]
\end{proof}

The weights in Theorem~\ref{thm:barrier} may have infinite integrals
both near the origin and at infinity. We therefore introduce the following relative functional.
Set
\begin{equation}\label{eq:relative-energy}
 \begin{aligned}
 \Phi_1(s)&=\int_1^s\frac{\dd a}{\phi(a)},\\
 \mathcal J(t)&=\int_{\Ct}f\Phi_1(h)\,\dd\sigma
       -\int_{\Sp}\int_1^{\rho(\omega,t)}G(a\omega)a^n\,\dd a\,\dd\omega.
 \end{aligned}
\end{equation}
By Lemma~\ref{lem:barriers}, there exist constants $0<c<C$,
independent of $t$, such that $c\le h,\rho\le C$.
Set
\[
 a:=\min\{1,c\},\qquad b:=\max\{1,C\},
\]
and define
\[
 m_\phi=\min_{s\in[a,b]}\phi(s)>0,\qquad
 M_G=\max_{\substack{s\in[a,b]\\\omega\in\Sp}}G(s\omega)<\infty.
\]
Then
\[
\begin{aligned}
 |\mathcal J(t)|
 &\le \int_{\Ct}f|\Phi_1(h)|\,\dd\sigma
   +\int_{\Sp}
       \left|\int_1^{\rho(\omega,t)}
                   G(s\omega)s^n\,\dd s\right|\dd\omega\\
 &\le \frac{b-a}{m_\phi}\int_{\Ct}f\,\dd\sigma
   +\frac{|\Sp|M_G}{n+1}
       \bigl(b^{n+1}-a^{n+1}\bigr)<\infty.
\end{aligned}
\]
The right-hand side is independent of $t$. Moreover, continuity of $h$ and $\rho$ at $t=0$ gives continuity
of $\mathcal J$ at $t=0$.

\begin{lemma}\label{lem:relative}
Assume the hypotheses of Theorem~\ref{thm:barrier}.
Let $h$ be a positive smooth strictly convex solution
to~\eqref{eq:scalar} with $\eta\equiv1$ and $h(\cdot,0)=s_0l$
on $\Ct\times[0,T)$.
Let $\mathcal J(t)$ be the functional defined
in~\eqref{eq:relative-energy}.
Then
\begin{equation}\label{eq:relative-dissipation}
 \mathcal J'(t)
 =-\int_{\Ct}\frac{G(X)}{hK}h_t^2\,\dd\sigma\le0,
 \qquad 0<t<T.
\end{equation}
Moreover, there exists a constant $C>0$, independent of $t$,
such that
\[
 |\mathcal J(t)|\le C,\qquad 0\le t<T,
\]
and
\[
 \int_0^T\int_{\Ct}\frac{G(X)}{hK}h_t^2
       \,\dd\sigma\,\dd t<\infty.
\]
\end{lemma}

\begin{proof}
By~\eqref{eq:first-variation} and~\eqref{eq:scalar} with $\eta=1$,
\[
\begin{aligned}
 \mathcal J'(t)
 &=\int_{\Ct}\frac f{\phi(h)}h_t\,\dd\sigma
   -\int_{\Ct}\frac{G(X)}K h_t\,\dd\sigma\\
 &=\int_{\Ct}\frac{G(X)}{hK}
       \big((h-h_t)-h\big)h_t\,\dd\sigma\\
 &=-\int_{\Ct}\frac{G(X)}{hK}h_t^2\,\dd\sigma\le0.
\end{aligned}
\]
The compact bounds for $h$ and $\rho$ in Lemma~\ref{lem:barriers}
give $|\mathcal J(t)|\le C$. Integrating the identity yields
\[
 \int_0^t\int_{\Ct}\frac{G(X)}{hK}h_s^2\,\dd\sigma\,\dd s
 =\mathcal J(0)-\mathcal J(t)
 \le\mathcal J(0)+C,\qquad 0<t<T.
\]
When $t\to T$ yields the desired estimate.
\end{proof}

\section{Curvature estimates and long-time existence}\label{sec:curvature}

The lower-order bounds obtained above allow the same curvature
argument. Throughout this section, $h$ is the solution of the flow~\eqref{eq:scalar}, and
\begin{equation}\label{eq:curvature-input}
 0<m\le h\le M,\qquad |\nabla h|\le C,\qquad
 0<\eta_-\le\eta(t)\le\eta_+,
\end{equation}
where $m,M,C,\eta_-$ and $\eta_+$ are positive constants
independent of $t$.
Denote
\[
 r:=\frac{fK}{\phi(h)G(X)},\qquad F:=\eta hr=h-h_t,
 \qquad\sigma_1=\tr A,\qquad S=\tr(A^{-1}).
\]
All components below are taken in local orthonormal frames for $\sigma$.
We also denote
\[
 \gamma:=\gradR\log G(X),\qquad
 \gamma_i=\langle\gamma,e_i\rangle,\qquad
 \gamma\cdot\nabla h=\sum_{i=1}^n\gamma_i h_i.
\]
All constants below are independent of $t$. Since $h_t=h(1-\eta r)$ and $h_\mu=\cot\theta\,h$, then we have
\[
 0=(h_t)_\mu-\cot\theta\,h_t
   =\partial_\mu\bigl(h(1-\eta r)\bigr)
       -\cot\theta\,h(1-\eta r)
   =-\eta h r_\mu.
\]
So
\begin{equation}\label{eq:r-boundary}
 r_\mu=0\qquad\text{on }\partial\Ct\times(0,T).
\end{equation}
This identity, together with the use of $r$ rather than $\eta r$ in
the auxiliary functions, avoids time derivatives of $\eta$ in the
curvature estimates.

\subsection{Bounds for the Gauss curvature}
\begin{lemma}\label{lem:K-upper}
Let $h$ be a positive smooth strictly convex solution
to~\eqref{eq:scalar} on $\Ct\times[0,T)$, with initial
support function $h_0$ and a positive time-dependent factor
$\eta(t)$.
Then the Gauss curvature $K$ satisfies
\[
 K(\xi,t)\le C,
 \qquad (\xi,t)\in\Ct\times[0,T),
\]
where $C>0$ is independent of $t$ and depends only on
$n,\theta,h_0,f,\phi,G$ and the constants in~\eqref{eq:curvature-input}.
\end{lemma}
\begin{proof}
Following the maximum-principle argument in
\cite[Lemma 4.3]{HHI}, set
\[
 \varepsilon=\frac m2,\qquad
 W=\frac{hr}{h-\varepsilon}
   =\frac{fhK}{(h-\varepsilon)\phi(h)G(X)}.
\]
The bounds~\eqref{eq:curvature-input} imply
\[
 h-\varepsilon\ge\frac m2>0,\qquad
 cK\le W\le CK.
\]

Fix $T'\in(0,T)$, and let $(\xi_0,t_0)$ be a maximum point
of $W$ on $\Ct\times[0,T']$.
If $t_0=0$, the maximum is bounded by the initial data.
We therefore consider $t_0>0$.

On $\partial\Ct\times(0,T)$, the boundary conditions
\eqref{eq:boundary} and~\eqref{eq:r-boundary} give
$h_\mu=\cot\theta\,h$ and $r_\mu=0$.
Hence
\begin{equation}\label{eq:W-boundary}
\begin{aligned}
 \frac{W_\mu}{W}
 &=\frac{h_\mu}{h}+\frac{r_\mu}{r}
       -\frac{h_\mu}{h-\varepsilon}\\
 &=\cot\theta-\frac{h\cot\theta}{h-\varepsilon}
 =-\frac{\varepsilon\cot\theta}{h-\varepsilon}<0.
\end{aligned}
\end{equation}
Thus, the maximum value of $W(\cdot, t_0)$ is attained at some point  $\xi_0\in\Ct\setminus\partial\Ct$.

Choose a local orthonormal frame $\{e_i\}_{i=1}^n$ for $\sigma$
such that $(A_{ij})$ is diagonal at $\xi_0$.
Then at $(\xi_0,t_0)$,
we have
\[
 0=\nabla_iW
 =\frac{(hr)_i}{h-\varepsilon}
  -\frac{hrh_i}{(h-\varepsilon)^2},
 \qquad
 (hr)_i=\frac{hr}{h-\varepsilon}h_i=Wh_i,
\]
and
\[
\begin{aligned}
 0\ge\nabla^2_{ii}W
 &=\frac{(hr)_{ii}}{h-\varepsilon}
   -\frac{2(hr)_ih_i+hrh_{ii}}{(h-\varepsilon)^2}
   +\frac{2hrh_i^2}{(h-\varepsilon)^3}\\
 &=\frac{(hr)_{ii}}{h-\varepsilon}
   -\frac{hrh_{ii}}{(h-\varepsilon)^2}
   -\frac{2hrh_i^2}{(h-\varepsilon)^3}
   +\frac{2hrh_i^2}{(h-\varepsilon)^3}\\
 &=\frac{(hr)_{ii}}{h-\varepsilon}
   -\frac{hrh_{ii}}{(h-\varepsilon)^2}.
\end{aligned}
\]
Since $h-\varepsilon>0$, it follows that
\[
 (hr)_{ii}
 \le\frac{hr}{h-\varepsilon}h_{ii}=Wh_{ii},
\]
hence
\[
\begin{aligned}
 (hr)_{ii}+hr
 &\le Wh_{ii}+(h-\varepsilon)W\\
 &=W(h_{ii}+h-\varepsilon)
 =W(A_{ii}-\varepsilon).
\end{aligned}
\]

Using $h_t=h-\eta hr$ and $K=(\det A)^{-1}$, we have
\[
\begin{aligned}
 K_t
 &=-K\sum_{i,j=1}^n A^{ij}(A_{ij})_t=-K\sum_{i=1}^n A^{ii}\bigl((h_t)_{ii}+h_t\bigr)\\
 &=-K\sum_{i=1}^n A^{ii}
       \bigl[A_{ii}-\eta\bigl((hr)_{ii}+hr\bigr)\bigr]\\
 &\le -K\sum_{i=1}^n A^{ii}
       \bigl[A_{ii}-\eta W(A_{ii}-\varepsilon)\bigr]\\
 &=K\left[
       -n+n\eta W
       -\eta\varepsilon W\sum_{i=1}^n A^{ii}
     \right],
\end{aligned}
\]
where we used
\[
 \sum_{i=1}^n A^{ii}A_{ii}=n.
\]
Moreover, the arithmetic-geometric mean inequality gives
\[
 \sum_{i=1}^n A^{ii}
 \ge n\left(\prod_{i=1}^n A^{ii}\right)^{1/n}
 =nK^{1/n}.
\]
Therefore,
\[
 \frac{K_t}{K}
 \le -n+n\eta W-n\eta\varepsilon WK^{1/n}.
\]

Next, since $f$ is independent of $t$, the evolution equation of
$W=hr/(h-\varepsilon)$ yields
\[
\begin{aligned}
 W_t
 &=\frac{(hr)_t}{h-\varepsilon}
   -\frac{hrh_t}{(h-\varepsilon)^2}\\
 &=W\left[
      \frac{K_t}{K}
      +\frac{h_t}{h}
      -\frac{\phi'(h)}{\phi(h)}h_t
      -\gradR\log G(X)\cdot X_t
      -\frac{h_t}{h-\varepsilon}
    \right]\\
 &=W\left[
      \frac{K_t}{K}
      -\left(
         \frac{\varepsilon}{h(h-\varepsilon)}
         +\frac{\phi'(h)}{\phi(h)}
       \right)h_t
      -\gradR\log G(X)\cdot X_t
    \right].
\end{aligned}
\]
We also have
\[
 (h_t)_i
 =h_i-\eta(hr)_i
 =(1-\eta W)h_i.
\]
Thus
\[
\begin{aligned}
 h_t&=h-\eta(h-\varepsilon)W,\\
 X_t&=\nabla h_t+h_t\nu\\
    &=(1-\eta W)\nabla h
      +\bigl(h-\eta(h-\varepsilon)W\bigr)\nu\\
    &=(1-\eta W)X+\eta\varepsilon W\nu.
\end{aligned}
\]
Using~\eqref{eq:curvature-input},
\[
 |h_t|+|X_t|\le C(1+W).
\]
The smoothness and positivity of $\phi$ and $G$ on the
corresponding compact sets also give
\[
 \left|
   \frac{\varepsilon}{h(h-\varepsilon)}
   +\frac{\phi'(h)}{\phi(h)}
 \right|
 +|\gradR\log G(X)|\le C.
\]
Substituting these estimates into the formula for $W_t$,
we obtain
\[
\begin{aligned}
 W_t
 &\le W\left(-n+n\eta W
             -n\eta\varepsilon WK^{1/n}\right)
       +CW(1+W)\\
 &\le C(W+W^2)
       -n\eta_-\varepsilon W^2K^{1/n}.
\end{aligned}
\]
Since $W\le CK$, we have $K^{1/n}\ge cW^{1/n}$.
At this positive-time maximum, $W_t(\xi_0,t_0)\ge0$.
Thus, for a constant $C\ge1$ independent of $T'$,
\begin{equation}\label{eq:W-reaction}
 0\le W_t\le C(W+W^2)-C^{-1}W^{2+1/n}.
\end{equation}
If $W(\xi_0,t_0)\ge1$, then
\[
 0\le2CW^2-C^{-1}W^{2+1/n},\qquad
 W(\xi_0,t_0)\le(2C^2)^n.
\]
Including the case $t_0=0$, we conclude that
\[
 \max_{\Ct\times[0,T']}W
 \le\max\left\{\max_{\Ct}W(\cdot,0),\,1,\,(2C^2)^n\right\}.
\]
Finally, $cK\le W$ and any $T'<T$ imply
\[
 K(\xi,t)\le C,
 \qquad (\xi,t)\in\Ct\times[0,T).
\]
This completes the proof.
\end{proof}

\begin{lemma}\label{lem:K-lower}
Let $h$ be a positive smooth strictly convex solution
to~\eqref{eq:scalar} on $\Ct\times[0,T)$, with initial support
function $h_0$.
Then we have
\[
 K(\xi,t)\ge c,
 \qquad (\xi,t)\in\Ct\times[0,T),
\]
where $c>0$ is independent of $t$ and depends only on
$n,\theta,h_0,f,\phi,G$ and the constants in~\eqref{eq:curvature-input}.
\end{lemma}
\begin{proof}
Denote
\[
 Q=-\log r-N_0\log h,\qquad
 r=\frac{fK}{\phi(h)G(X)},
\]
where $N_0>1$ will be chosen later.
Fix $T'\in(0,T)$, and let $(\xi_0,t_0)$ be a maximum point
of $Q$ on $\Ct\times[0,T']$.
If $t_0=0$, the maximum is bounded by the initial data.
We therefore consider $t_0>0$.
On $\partial\Ct\times(0,T)$,
\eqref{eq:boundary} and~\eqref{eq:r-boundary} give
\[
 Q_\mu=-\frac{r_\mu}{r}-N_0\frac{h_\mu}{h}
       =-N_0\cot\theta<0.
\]
Thus, the maximum value of $Q(\cdot, t_0)$ is attained at some point  $\xi_0\in\Ct\setminus\partial\Ct$.
Moreover,
\[
 Q_t(\xi_0,t_0)
 =\lim_{\delta\to 0}
   \frac{Q(\xi_0,t_0)-Q(\xi_0,t_0-\delta)}{\delta}
 \ge0.
\]

Choose a local orthonormal frame for $\sigma$ such that
$(A_{ij})$ is diagonal at $\xi_0$.
Then at $(\xi_0,t_0)$, we have
\[
 0=\nabla_iQ=-\frac{r_i}{r}-N_0\frac{h_i}{h},
\]
so
\[
 r_i=-N_0r\frac{h_i}{h}.
\]
We also have
\[
\begin{aligned}
 0\ge\nabla^2_{ii}Q
 &=-\frac{r_{ii}}r+\frac{r_i^2}{r^2}
   -N_0\frac{h_{ii}}h+N_0\frac{h_i^2}{h^2}\\
 &=-\frac{r_{ii}}r-N_0\frac{h_{ii}}h
   +N_0(N_0+1)\frac{h_i^2}{h^2}.
\end{aligned}
\]
Hence
\[
 -r_{ii}
 \le N_0r\frac{h_{ii}}h
      -N_0(N_0+1)r\frac{h_i^2}{h^2}.
\]
Using $h_t=h-\eta hr$, we obtain that
\[
\begin{aligned}
 (h_t)_{ii}
 &=(1-\eta r)h_{ii}-2\eta h_ir_i-\eta h r_{ii}\\
 &\le(1-\eta r)h_{ii}
      +2N_0\eta r\frac{h_i^2}{h}
      +N_0\eta r h_{ii}
      -N_0(N_0+1)\eta r\frac{h_i^2}{h}\\
 &=\bigl[1+(N_0-1)\eta r\bigr]h_{ii}
      -N_0(N_0-1)\eta r\frac{h_i^2}{h}.
\end{aligned}
\]
Since $h_{ii}=A_{ii}-h$ and $h_t=h(1-\eta r)$,
\[
\begin{aligned}
 (h_t)_{ii}+h_t
 &\le\bigl[1+(N_0-1)\eta r\bigr](A_{ii}-h)
       +h(1-\eta r)
       -N_0(N_0-1)\eta r\frac{h_i^2}{h}\\
 &=\bigl[1+(N_0-1)\eta r\bigr]A_{ii}
       -N_0\eta hr
       -N_0(N_0-1)\eta r\frac{h_i^2}{h}.
\end{aligned}
\]
Multiplying by $A^{ii}$ and summing over $i$ yields
\[
\begin{aligned}
 \sum_{i=1}^n A^{ii}\bigl((h_t)_{ii}+h_t\bigr)
 &\le n+n(N_0-1)\eta r
       -N_0\eta hr\sum_{i=1}^n A^{ii}\\
 &\quad-\frac{N_0(N_0-1)\eta r}{h}
                 \sum_{i=1}^n A^{ii}h_i^2.
\end{aligned}
\]

On the other hand, we know that
\[
 (h_t)_i=(1-\eta r)h_i-\eta hr_i
        =\bigl[1+(N_0-1)\eta r\bigr]h_i.
\]
Therefore,
\[
\begin{aligned}
 X_t
 &=\nabla h_t+h_t\nu\\
 &=\bigl[1+(N_0-1)\eta r\bigr]\nabla h
       +h(1-\eta r)\nu\\
 &=(1-\eta r)X+N_0\eta r\nabla h.
\end{aligned}
\]
Writing $\gamma=\gradR\log G(X)$, we calculate the following evolution equation
\[
\begin{aligned}
 Q_t
 &=-\frac{r_t}{r}-N_0\frac{h_t}{h}\\
 &=-\frac{K_t}{K}
    +\frac{\phi'(h)}{\phi(h)}h_t
    +\gamma\cdot X_t-N_0\frac{h_t}{h}\\
 &=\sum_{i=1}^n A^{ii}\bigl((h_t)_{ii}+h_t\bigr)
    +\frac{\phi'(h)}{\phi(h)}h_t
    +\gamma\cdot X_t-N_0\frac{h_t}{h}.
\end{aligned}
\]
Substituting the preceding formulas gives
\[
\begin{aligned}
 Q_t
 &\le n+n(N_0-1)\eta r
       -N_0\eta hr\sum_{i=1}^n A^{ii}
       -\frac{N_0(N_0-1)\eta r}{h}
                    \sum_{i=1}^n A^{ii}h_i^2\\
 &\quad+\left(
          \frac{h\phi'(h)}{\phi(h)}+\gamma\cdot X-N_0
        \right)(1-\eta r)
       +N_0\eta r\,\gamma\cdot\nabla h.
\end{aligned}
\]
Recalling $E(h,X)=n+h\phi'(h)/\phi(h)+\gamma\cdot X$, we obtain
\begin{equation}\label{eq:K-lower-reaction}
\begin{aligned}
 0\le Q_t
 &\le E(h,X)-N_0\\
 &\quad+\eta r\left[
       N_0(n+1)-E(h,X)+N_0\gamma\cdot\nabla h
       -N_0h\sum_{i=1}^n A^{ii}\right.\\
 &\hspace{37mm}\left.
       -\frac{N_0(N_0-1)}h\sum_{i=1}^n A^{ii}h_i^2
       \right].
\end{aligned}
\end{equation}

By~\eqref{eq:curvature-input}, we have
\[
 m\le |X|\le R:=\sqrt{M^2+C^2}.
\]
Choose $N_0$ such that
\[
 N_0>1+
 \max_{\substack{m\le s\le M,\ Y\in\Hh\\m\le |Y|\le R}}
       |E(s,Y)|.
\]

Then $E(h,X)-N_0\le-1$.
Since $A^{ii}>0$ and $N_0>1$, the last two terms in brackets
in~\eqref{eq:K-lower-reaction} are nonpositive.
Then we have
\[
 0\le Q_t\le-1+Cr
 =-1+Ce^{-Q}h^{-N_0}
 \le-1+Cm^{-N_0}e^{-Q},
\]
where $C>0$ is independent of $t$.
Thus
\[
 Q(\xi,t)\leq Q(\xi_0,t_0)\le\ln(C m^{-N_0}).
\]
Including the case $t_0=0$ and any $T'<T$,
we obtain
\[
 Q\le
 \max\left\{
       \max_{\Ct}Q(\cdot,0),\,
       \ln(Cm^{-N_0}),\,1
     \right\}=:\widetilde{C}.
\]
So,
\[
 r=e^{-Q}h^{-N_0}\ge e^{-\widetilde{C}}M^{-N_0},
\]
and therefore
\[
 K=\frac{r\phi(h)G(X)}f
 \ge e^{-\widetilde{C}}M^{-N_0}
 \frac{\displaystyle
       \min_{s\in[m,M]}\phi(s)\,
       \cdot\min_{\substack{Y\in\Hh\\m\le |Y|\le R}}G(Y)}
      {\displaystyle\max_{\Ct}f}
 =:c>0.
\]

This completes the proof.
\end{proof}
\begin{remark}
The proof uses boundedness of $E(h,X)$ on the compact range of the
flow; its negativity in~\eqref{ass:U} is not required here.
\end{remark}

In particular, the preceding lemmas give
\begin{equation}\label{eq:det-speed-bounds}
 c\le\det A\le C,\qquad c\le F\le C.
\end{equation}
The positive lower bound for $F$ will be used in the trace estimate.

\subsection{The principal radii}

We retain the notation $A,F,\sigma_1,S$ introduced above.
At a diagonal point, let $\lambda_i>0$ be the eigenvalues of $A$,
and $\kappa_i=\lambda_i^{-1}$ be the principal curvatures.
We write $\sigma_k(A)$ for the $k$-th elementary symmetric
functions, so that $\sigma_{n-1}=(\det A)S$.
Set
\[
 \psi=\log F+\log\det A
 =\log\eta(t)+\log f+\log h-\log\phi(h)-\log G(X),
\]
\[
 \mathcal L=A^{ij}\nabla_{ij},\qquad
 \mathcal P=\partial_t-F\mathcal L.
\]
At a diagonal point, define
\begin{equation}\label{eq:W-weights}
 W_i=\frac1{\sigma_1}+\frac1{\lambda_i}
                      -\frac1{S\lambda_i^2}.
\end{equation}
All sums run from $1$ to $n$ unless otherwise indicated.
Throughout the proof, $C$ denotes a positive constant independent
of $t$, which may change from line to line.

\begin{proposition}\label{prop:curvature}\label{lem:symmetric-evolution}
Let $h$ be a positive smooth strictly convex solution
to~\eqref{eq:scalar} on $\Ct\times[0,T)$, with initial support
function $h_0$.
Then we have
\begin{equation}\label{eq:curvature-bounds}
 C^{-1}\sigma\le\nabla^2h+h\sigma\le C\sigma
 \qquad\text{on }\Ct\times[0,T),
\end{equation}
where $C>1$ depends only on $n,\theta,h_0,f,\phi,G$
and the constants in~\eqref{eq:curvature-input}.
In particular,
\[
 C^{-1}\le\lambda_i,\kappa_i\le C,\qquad
 |\nabla^2h|\le C.
\]
\end{proposition}

\begin{proof}
By Lemma~\ref{lem:K-upper} and Lemma~\ref{lem:K-lower},
\[
 C^{-1}\le\det A,F\le C,\qquad
 \sigma_1\ge n(\det A)^{1/n}\ge C^{-1}.
\]
Following~\cite[Section 4]{Gao} and~\cite[Lemma 4.5]{HHI}, consider the auxiliary function
\begin{equation}\label{eq:trace-test}
 Q=\log\sigma_1+\log\sigma_{n-1}+ah+b|\nabla h|^2,
 \qquad a<0,\quad b>0.
\end{equation}
Fix $T'\in(0,T)$, and let $(\xi_0,t_0)$ maximize $Q$ on
$\Ct\times[0,T']$. The case $t_0=0$ is controlled by $h_0$.
For $t_0>0$, the maximum condition gives
\[
 Q_t(\xi_0,t_0)
 =\lim_{\delta\downarrow0}
   \frac{Q(\xi_0,t_0)-Q(\xi_0,t_0-\delta)}{\delta}\ge0.
\]

\smallskip\noindent\emph{Case 1: $\xi_0\in\Ct\setminus\partial\Ct$.}
Choose a local orthonormal frame diagonalizing $A$ at $\xi_0$.
All calculations in this case are at $(\xi_0,t_0)$.
Since $\log\sigma_{n-1}=\log\det A+\log S$, we have
\begin{equation}\label{eq:Q-stationary}
\begin{aligned}
 0=Q_k
 &=\frac{\sum_iA_{ii,k}}{\sigma_1}
   +\sum_i\frac{A_{ii,k}}{\lambda_i}
   -\frac1S\sum_i\frac{A_{ii,k}}{\lambda_i^2}
   +\bigl(a+2b(\lambda_k-h)\bigr)h_k\\
 &=\sum_iW_iA_{ii,k}
   +\bigl(a+2b(\lambda_k-h)\bigr)h_k.
\end{aligned}
\end{equation}
For each fixed $k$,
\[
\begin{aligned}
 0\ge Q_{kk}
 &=\frac{(\sigma_1)_{kk}}{\sigma_1}
   -\frac{((\sigma_1)_k)^2}{\sigma_1^2}
   +(\log\det A)_{kk}+\frac{S_{kk}}S-\frac{S_k^2}{S^2}
   +ah_{kk}+b(|\nabla h|^2)_{kk}.
\end{aligned}
\]
Thus $\mathcal LQ\le0$ and $\mathcal PQ\ge0$.

Since $S\lambda_i\ge1$, direct summation gives
\begin{equation}\label{eq:W-sums}
 W_i\ge0,\qquad
 \sum_iW_i\lambda_i=n,\qquad
 \sum_iW_i\le\frac n{\sigma_1}+S,\qquad
 \sum_iW_i\lambda_i^2\le2\sigma_1.
\end{equation}
Differentiating $\log\det A$ and $S$, for each fixed $k$ we have
\[
 (\log\det A)_k=\sum_i\frac{A_{ii,k}}{\lambda_i},\qquad
 S_k=-\sum_i\frac{A_{ii,k}}{\lambda_i^2},
\]
\[
\begin{aligned}
 (\log\det A)_{kk}
 &=\sum_i\frac{A_{ii,kk}}{\lambda_i}
       -\sum_{i,j}\frac{A_{ij,k}^2}{\lambda_i\lambda_j},\\
 S_{kk}
 &=-\sum_i\frac{A_{ii,kk}}{\lambda_i^2}
       +2\sum_{i,j}\frac{A_{ij,k}^2}{\lambda_i^2\lambda_j}.
\end{aligned}
\]
Using $A_{ii,kk}=A_{kk,ii}+\lambda_i-\lambda_k$, we obtain
\[
\begin{aligned}
 \mathcal L\sigma_1
 &=\Delta\log\det A
   +\sum_{i,j,k}\frac{A_{ij,k}^2}{\lambda_i\lambda_j}
   +\sigma_1S-n^2,\\
 \mathcal LS
 &=-\sum_i\frac{(\log\det A)_{ii}}{\lambda_i^2}
   +\sum_{i,j,k}\frac{A_{ij,k}^2}{\lambda_i\lambda_j^2\lambda_k}
   -S^2+n\sum_i\lambda_i^{-2}.
\end{aligned}
\]
By Cauchy--Schwarz inequality and the full symmetry of $A_{ij,k}$,
\[
\begin{aligned}
 \sum_k\frac{((\sigma_1)_k)^2}{\lambda_k}
 &\le\sigma_1\sum_{i,k}\frac{A_{ii,k}^2}{\lambda_i\lambda_k}
 \le\sigma_1\sum_{i,j,k}\frac{A_{ij,k}^2}{\lambda_i\lambda_j},\\
 \sum_k\frac{S_k^2}{\lambda_k}
 &\le S\sum_{i,k}\frac{A_{ii,k}^2}{\lambda_i^3\lambda_k}
 \le S\sum_{i,j,k}\frac{A_{ij,k}^2}{\lambda_i\lambda_j^2\lambda_k}.
\end{aligned}
\]
Applying these inequalities to
$\mathcal L(\log\sigma_1+\log\det A+\log S)$ yields
\begin{equation}\label{eq:spatial-B}
 \mathcal L\log(\sigma_1\cdot \sigma_{n-1})
 \ge\sum_iW_i(\log\det A)_{ii}
          -\frac{n^2}{\sigma_1}+\frac nS\sum_i\lambda_i^{-2}.
\end{equation}
On the other hand, $h_t=h-F$ gives
\[
\begin{aligned}
 (A_{ij})_t&=A_{ij}-F_{ij}-F\sigma_{ij},\\
 \partial_t\log(\sigma_1\cdot\sigma_{n-1})
 &=n-\sum_iW_iF_{ii}-F\sum_iW_i,\\
 \frac{F_{ii}}F
 &=\psi_{ii}-(\log\det A)_{ii}+((\log F)_i)^2.
\end{aligned}
\]
Together with~\eqref{eq:spatial-B}, this implies
\begin{equation}\label{eq:symmetric-evolution}
\begin{aligned}
 \frac{\mathcal P\log(\sigma_1\sigma_{n-1})}F
 &\le\frac nF+\frac{n(n-1)}{\sigma_1}-S
       -\frac{n-1}S\sum_i\lambda_i^{-2}\\
 &\quad-\sum_iW_i\psi_{ii}-\sum_iW_i((\log F)_i)^2.
\end{aligned}
\end{equation}

For the gradient term, the Codazzi equation gives
$h_{jk,k}=A_{kk,j}-h_k\delta_{jk}$, and hence
\[
\begin{aligned}
 \mathcal L|\nabla h|^2
 &=2\sum_k\frac{(\lambda_k-h)^2}{\lambda_k}
    +2\langle\nabla h,\nabla\log\det A\rangle
    -2\sum_k\frac{h_k^2}{\lambda_k},\\
 \frac{\partial_t|\nabla h|^2}F
 &=\frac{2|\nabla h|^2}F
    -2\langle\nabla h,\nabla\psi\rangle
    +2\langle\nabla h,\nabla\log\det A\rangle.
\end{aligned}
\]
Subtracting these identities gives
\begin{equation}\label{eq:h-gradient-evolution}
 \frac{\mathcal P|\nabla h|^2}F
 =\frac{2|\nabla h|^2}F-2\langle\nabla h,\nabla\psi\rangle
       -2\sigma_1+4nh-2h^2S
       +2\sum_i\frac{h_i^2}{\lambda_i}.
\end{equation}
Also,
\begin{equation}\label{eq:h-evolution}
 \frac{\mathcal Ph}F
 =\frac hF-1-\sum_i\frac{\lambda_i-h}{\lambda_i}
 =\frac hF-(n+1)+hS.
\end{equation}

By~\eqref{eq:composition},
$D_iX=\lambda_i e_i$ and
$D^2_{ii}X=\sum_kA_{ii,k}e_k-\lambda_i\nu$.
Since $\eta$ depends only on $t$, for each fixed $i$ we obtain
\begin{equation}\label{eq:psi-derivatives}
\begin{aligned}
 \psi_i
 &=(\log f)_i
   +\left(\frac1h-\frac{\phi'(h)}{\phi(h)}\right)h_i
   -\lambda_i\gamma_i,\\
 \psi_{ii}
 &=(\log f)_{ii}
   +\left[-\frac1{h^2}-\frac{\phi''(h)}{\phi(h)}
                  +\frac{\phi'(h)^2}{\phi(h)^2}\right]h_i^2\\
 &\quad+\left(\frac1h-\frac{\phi'(h)}{\phi(h)}\right)(\lambda_i-h)
   -\lambda_i^2(\hessR\log G)(X)(e_i,e_i)\\
 &\quad-\sum_k\gamma_kA_{ii,k}
       +\lambda_i\langle\gamma,\nu\rangle.
\end{aligned}
\end{equation}
Using~\eqref{eq:Q-stationary}, we have
\begin{equation}\label{eq:G-cancellation}
\begin{aligned}
 &\sum_k\gamma_k\sum_iW_iA_{ii,k}
       +2b\sum_k\lambda_kh_k\gamma_k\\
 &=-\sum_k\bigl(a+2b(\lambda_k-h)\bigr)h_k\gamma_k
       +2b\sum_k\lambda_kh_k\gamma_k\\
 &=(2bh-a)\langle\gamma,\nabla h\rangle.
\end{aligned}
\end{equation}
The remaining terms are bounded using~\eqref{eq:W-sums}:
\[
\begin{aligned}
 \sum_iW_i(1+h_i^2+|\lambda_i-h|+\lambda_i)
 &\le C(1+S),\\
 \left|\sum_iW_i\lambda_i^2(\hessR\log G)(X)(e_i,e_i)\right|
 &\le C\sigma_1,\\
 |(2bh-a)\langle\gamma,\nabla h\rangle|
 &\le C(|a|+b),\\
 2b\left|\left\langle\nabla h,
       \nabla\log f+\left(\frac1h-\frac{\phi'(h)}{\phi(h)}\right)
                    \nabla h\right\rangle\right|
 &\le Cb.
\end{aligned}
\]
Therefore
\begin{equation}\label{eq:weighted-psi-estimate}
 -\sum_iW_i\psi_{ii}-2b\langle\nabla h,\nabla\psi\rangle
 \le C(1+|a|+b)+C\sigma_1+CS,
\end{equation}
where $C$ is independent of $a,b,T'$.
Combining~\eqref{eq:symmetric-evolution}--\eqref{eq:h-evolution}
and~\eqref{eq:weighted-psi-estimate}, and dropping the nonpositive
square terms, we obtain
\begin{equation}\label{eq:trace-reaction}
\begin{aligned}
 0\le\frac{\mathcal PQ}F
 &\le C(1+|a|+b)+(C-2b)\sigma_1
               +(ah-2bh^2+2bC^2+C-1)S\\
 &\le C(1+|a|+b)+(C-2b)\sigma_1
               +(am+2bC^2+C)S.
\end{aligned}
\end{equation}
Here we used $\sum_i h_i^2/\lambda_i\le C^2S$.
Fix the $C$ in~\eqref{eq:trace-reaction} large enough,
choose $b>0$ and then $a<0$ such that
\[
 2b-C\ge1,\qquad am+2bC^2+C\le-1.
\]
Thus $\sigma_1+S\le C(1+|a|+b)$ at an interior maximum.
We keep $b$ fixed and fix $a$ after the boundary estimate;
decreasing $a$ preserves these inequalities.

\smallskip\noindent\emph{Case 2: $\xi_0\in\partial\Ct$.}
Take $e_n=\mu$ and diagonalize the tangential block of $A$.
By~\eqref{eq:boundary}, $A_{\alpha n}=h_{\alpha n}=0$,
so the whole matrix is diagonal and $\lambda_n=A_{nn}$.
At this maximum point, $Q_\mu\ge0$.

Using $\nabla_{e_\alpha}\mu=\cot\theta e_\alpha$,
the tangential derivative of $A(e_\alpha,\mu)=0$
and the Codazzi equation give
\[
 A_{\alpha\alpha,\mu}=A_{\alpha n,\alpha}
 =\cot\theta(\lambda_n-\lambda_\alpha).
\]
Since $r_\mu=0$, $h_\mu=\cot\theta h$ and
$D_\mu X=\lambda_n\mu$,
\begin{equation}\label{eq:boundary-determinant}
\begin{aligned}
 (\log\det A)_\mu
 &=(\log f)_\mu-\cot\theta\frac{h\phi'(h)}{\phi(h)}
       -\lambda_n\langle\gamma,\mu\rangle,\\
 |(\log\det A)_\mu|&\le C(1+\lambda_n).
\end{aligned}
\end{equation}
Moreover,
\[
\begin{aligned}
 A_{nn,\mu}
 &=\lambda_n\left(
      (\log\det A)_\mu
      -\sum_{\alpha<n}\frac{A_{\alpha\alpha,\mu}}{\lambda_\alpha}
    \right)\\
 &=\lambda_n(\log\det A)_\mu
      -\cot\theta\lambda_n(\lambda_nS-n).
\end{aligned}
\]
Consequently,
\[
\begin{aligned}
 (\sigma_1)_\mu
 &=\lambda_n(\log\det A)_\mu
       +\cot\theta(2n\lambda_n-\sigma_1-\lambda_n^2S),\\
 S_\mu
 &=-\sum_{\alpha<n}
       \frac{\cot\theta(\lambda_n-\lambda_\alpha)}{\lambda_\alpha^2}
       -\frac{A_{nn,\mu}}{\lambda_n^2}\\
 &=-\frac{(\log\det A)_\mu}{\lambda_n}
       -\cot\theta\lambda_n
           \sum_{\alpha<n}(\lambda_\alpha^{-1}-\lambda_n^{-1})^2,\\
 (|\nabla h|^2)_\mu
 &=2h_nh_{nn}
  =2\cot\theta h(\lambda_n-h).
\end{aligned}
\]
Using $\sigma_{n-1}=(\det A)S$, we obtain
\begin{equation}\label{eq:boundary-sigma}
 (\log\sigma_{n-1})_\mu
 =\left(1-\frac1{\lambda_nS}\right)(\log\det A)_\mu
   -\frac{\cot\theta\lambda_n}S
      \sum_{\alpha<n}(\lambda_\alpha^{-1}-\lambda_n^{-1})^2.
\end{equation}
Substitution into $Q_\mu\ge0$ gives
\begin{equation}\label{eq:boundary-Q}
\begin{aligned}
 0\le\frac{\sigma_1}{\lambda_n}Q_\mu
 ={}&\left(1+\frac{\sigma_1}{\lambda_n}
               -\frac{\sigma_1}{\lambda_n^2S}\right)
          (\log\det A)_\mu\\
 &+\cot\theta\left(
       2n-\frac{\sigma_1}{\lambda_n}-\lambda_nS
     \right)
       +a\cot\theta h\frac{\sigma_1}{\lambda_n}\\
 &+2b\cot\theta h(\lambda_n-h)\frac{\sigma_1}{\lambda_n}
       -\frac{\cot\theta\sigma_1}S
          \sum_{\alpha<n}(\lambda_\alpha^{-1}-\lambda_n^{-1})^2.
\end{aligned}
\end{equation}
Here
\[
\begin{aligned}
 1&\le1+\frac{\sigma_1}{\lambda_n}
              -\frac{\sigma_1}{\lambda_n^2S}
       \le1+\frac{\sigma_1}{\lambda_n},\\
 \frac{\sigma_1}{\lambda_n}+\lambda_nS
 &=\sum_i\left(
       \frac{\lambda_i}{\lambda_n}+\frac{\lambda_n}{\lambda_i}
     \right)\ge2n.
\end{aligned}
\]

If $\lambda_n\le m/2$, then $\lambda_n-h<0$ and
$|(\log\det A)_\mu|\le C$.
Further decrease $a$ so that $am\cot\theta\le-2C-1$.
Equation~\eqref{eq:boundary-Q} then gives
\[
 0\le C+(C+am\cot\theta)\frac{\sigma_1}{\lambda_n}
   \le C-(C+1)\frac{\sigma_1}{\lambda_n}\le-1,
\]
a contradiction.

If $\lambda_n>m/2$, then
\[
\begin{aligned}
 \left(1+\frac{\sigma_1}{\lambda_n}\right)(1+\lambda_n)
 &\le C(1+\sigma_1),\\
 2b\cot\theta h(\lambda_n-h)\frac{\sigma_1}{\lambda_n}
 &\le2bM\cot\theta\sigma_1.
\end{aligned}
\]
Using these estimates and Cauchy--Schwarz inequality
in~\eqref{eq:boundary-Q}, we obtain
\[
\begin{aligned}
 0
 &\le C(1+b)(1+\sigma_1)
       -\frac{\cot\theta\sigma_1}S
          \sum_{\alpha<n}
             (\lambda_\alpha^{-1}-\lambda_n^{-1})^2\\
 &\le C(1+b)(1+\sigma_1)
       -\frac{\cot\theta}{n-1}\frac{\sigma_1}S
          \left(S-\frac n{\lambda_n}\right)^2.
\end{aligned}
\]
When $S\ge4n/m$, we have $n/\lambda_n<2n/m\le S/2$, so
\[
 0\le C(1+b)(1+\sigma_1)
       -\frac{\cot\theta}{4(n-1)}\sigma_1S,
 \qquad
 S\le C(1+b)(1+\sigma_1^{-1})\le C(1+b).
\]
The same bound holds when $S<4n/m$.
Since $b$ is fixed, $S\le C$. The determinant bounds now give
\[
 \lambda_i\ge S^{-1}\ge C^{-1},\qquad
 \lambda_i=\frac{\det A}{\prod_{j\ne i}\lambda_j}
           \le CS^{n-1}\le C.
\]

Fix one $a<0$ satisfying the interior and boundary requirements.
At the maximum point, $\sigma_1$ and
$\sigma_{n-1}\le n\sigma_1^{n-1}$ are bounded, hence $Q\le C$.
At every point of $\Ct\times[0,T']$,
\[
 \sigma_{n-1}=(\det A)S
 \ge n(\det A)^{(n-1)/n}\ge C^{-1},
\]
and therefore
\[
 \log\sigma_1
 =Q-ah-b|\nabla h|^2-\log\sigma_{n-1}
 \le C+|a|M.
\]
Thus
\[
 \lambda_{\max}\le\sigma_1\le C,\qquad
 \lambda_{\min}
 \ge\frac{\det A}{\lambda_{\max}^{n-1}}\ge C^{-1}.
\]
Since $T'<T$ is arbitrary,~\eqref{eq:curvature-bounds} follows.
The estimates for $\kappa_i=\lambda_i^{-1}$ and
$\nabla^2h=A-h\sigma$ follow immediately.
\end{proof}

\begin{remark}
(a)~The estimate in Proposition~\ref{prop:curvature} can also be
obtained using the auxiliary function
\[
 \widetilde Q=\log\sigma_1(A)-ah+b|\nabla u|^2,
 \qquad u=\frac hl.
\]
Here $b>0$ is chosen sufficiently large, and then $a>0$ is
chosen sufficiently large.
We omit this alternative proof.

(b)~The boundary estimate uses $\cot\theta>0$ and does not cover
$\theta=\pi/2$ for general $G$.
For $n=2$, $\sigma_{n-1}=\sigma_1$, and the same proof applies.
\end{remark}

\subsection{Higher-order estimates}
\begin{proposition}\label{prop:longtime}
Under the assumptions of either Theorem~\ref{thm:even} or
Theorem~\ref{thm:barrier}, the solution exists for all $t>0$.
For every $\delta>0$ and integer $k\ge0$, its spatial $C^k$ norm
is uniformly bounded for $t\ge\delta$. 
\end{proposition}
\begin{proof}
Lemma~\ref{lem:even-lower} and Lemma~\ref{lem:barriers} give
\eqref{eq:curvature-input}. Proposition~\ref{prop:curvature} then implies
\[
 c\sigma\le A\le C\sigma,\qquad c\le F\le C,\qquad
 c|\zeta|^2\le FA^{ij}\zeta_i\zeta_j\le C|\zeta|^2.
\]
In particular, $|h_t|=|h-F|\le C$. When $\eta=N/D$,
\eqref{eq:eta-variation} with $v=h_t$ gives
\begin{equation}\label{eq:eta-time}
 \frac{\dd\eta}{\dd t}
 =\frac1D\int_{\Ct}
 \left[\frac{(n+1)G(X)+\gradR G(X)\cdot X}{K}
       -\eta\frac f{\phi(h)}
          \left(1-\frac{h\phi'(h)}{\phi(h)}\right)\right]
 h_t\,\dd\sigma.
\end{equation}
All coefficients are bounded and $D\ge c$, so
\[
 |\eta_t|\le C\int_{\Ct}|h_t|\,\dd\sigma\le C.
\]
For $\eta=1$, we have $\eta_t=0$.

Set $h=lu$. From~\eqref{eq:l-u}, \eqref{eq:boundary} and the preceding bounds,
\[
 u_\mu=0,\qquad
 u_t=u-\frac{\eta fu}{\phi(lu)G(l\nabla u+u\xi)\det A},\qquad
 \|u\|_{C^2}+\|u_t\|_{C^0}\le C.
\]
At fixed $\xi,t,u,\nabla u$, the matrix $A$ depends affinely on
$\nabla^2u$ with positive coefficient $l$. For any symmetric matrix $B$,
\[
\begin{aligned}
 \left.\frac{\dd^2}{\dd s^2}\right|_{s=0}
       \bigl(-\det(A+sB)^{-1}\bigr)
 =-\frac{(\tr(A^{-1}B))^2+\tr(A^{-1}BA^{-1}B)}{\det A}\le0.
\end{aligned}
\]
Thus the operator is concave in $\nabla^2u$ and uniformly parabolic.
Also, for $|t-s|\le1$, interpolation gives
\[
\begin{aligned}
 \|u(t)-u(s)\|_{C^0}&\le C|t-s|,\\
 \|\nabla u(t)-\nabla u(s)\|_{C^0}
 &\le C\|u(t)-u(s)\|_{C^0}^{1/2}
          \|u(t)-u(s)\|_{C^2}^{1/2}
       +C\|u(t)-u(s)\|_{C^0}\\
 &\le C|t-s|^{1/2},\qquad
 |\eta(t)-\eta(s)|\le C|t-s|.
\end{aligned}
\]
In fixed boundary-normal coordinates, the boundary condition is
$u_{x_n}=0$. Freezing $u$, $\nabla u$ and $\eta$ in the lower-order
coefficients gives a concave operator with a uniform parabolic
H\"older modulus in space and time. On the compact positive-definite
matrix range, the ellipticity constants are uniform.
The function $A\mapsto-(\det A)^{-1}$ admits a smooth concave
uniformly elliptic extension. Take its supporting-affine envelope
on a larger compact positive-definite matrix interval, mollify it,
and use a smooth cutoff to retain the original function near the
actual matrix range. In the gluing region the mollification converges
to the original function in $C^2$; strict concavity and the ellipticity
bounds therefore persist for a sufficiently small smoothing radius.
The interior estimates and the parabolic boundary estimate described
in~\cite[Theorem~5.4 and Section~7, p.~149]{Lieberman2002} therefore give,
for some $\alpha\in(0,1)$,
\[
 \|h\|_{C^{2+\alpha,1+\alpha/2}
 (\Ct\times[t,\min\{t+1,\tau\}])}
 \le C_\delta,\qquad \delta\le t<\tau<T.
\]
The estimates are applied on a larger interval starting at
$t-\delta/2$, so their constants are independent of $t$ and $\tau$.
For the linear oblique Schauder estimates used below, see
\cite{Lieberman,LiebermanErrata}.

For $\eta=N/D$, the integrand in~\eqref{eq:eta-time} is
$C^{\alpha,\alpha/2}$, so $\eta_t\in C^{\alpha/2}$ and
$\eta\in C^{1+\alpha/2}$. On successively smaller positive-time
intervals, differentiating the equation and the integral formula,
then applying oblique Schauder estimates, gives the induction
\[
\begin{aligned}
 h\in C^{2j+2+\alpha,j+1+\alpha/2}
 &\Longrightarrow \eta\in C^{j+1+\alpha/2}\\
 &\Longrightarrow h\in C^{2j+4+\alpha,j+2+\alpha/2},
 \qquad j=0,1,\ldots.
\end{aligned}
\]
For $\eta=1$ the time factor is constant. Consequently,
\[
 \|h\|_{C^{k+\alpha,(k+\alpha)/2}
 (\Ct\times[t,\min\{t+1,\tau\}])}
 \le C_{k,\delta},\qquad k\ge2,\quad\delta\le t<\tau<T.
\]
The constants are independent of $t$ and $\tau$.
If $T<\infty$, let $\tau\uparrow T$ in these estimates to obtain
\[
 \|h(t)-h(s)\|_{C^k}
 \le\int_s^t\|h_\tau\|_{C^k}\,\dd\tau\le C_k|t-s|,
\]
and hence
\[
 h(\cdot,t)\longrightarrow h_T\quad\text{in every }C^k(\Ct),
 \qquad h_T\ge c,\quad A[h_T]\ge c\sigma,\quad
 (h_T)_\mu=\cot\theta\,h_T.
\]
Lemma~\ref{lem:short-time} extends the solution past $T$, a contradiction.
Thus $T=\infty$, with the stated uniform estimates.
\end{proof}

\section{Convergence and proofs of the main results}\label{sec:convergence}

\subsection{Stationary subsequences}
We now complete the existence results. The two dissipation identities
have the same geometric integrand; only the choice of the monotone
functional and the scalar coefficient differ.

\begin{proof}[Proofs of Theorems~\ref{thm:even} and~\ref{thm:barrier}]
Proposition~\ref{prop:longtime} gives existence for all $t>0$ and
uniform positive-time estimates. Put
\[
 \mathcal D(t)=\int_{\Ct}\frac{G(X)}{hK}h_t^2\,\dd\sigma.
\]
By~\eqref{eq:entropy} and~\eqref{eq:relative-dissipation}, respectively,
\[
\begin{aligned}
 \int_0^R\mathcal D(t)\,\dd t
 &\le\eta_+\bigl(J_+(0)-J_+(R)\bigr)\le C
       &&\text{under }\eqref{ass:plus},\\
 \int_0^R\mathcal D(t)\,\dd t
 &\le\eta_+\bigl(J_-(R)-J_-(0)\bigr)\le C
       &&\text{under }\eqref{ass:minus},\\
 \int_0^R\mathcal D(t)\,\dd t
 &=\mathcal J(0)-\mathcal J(R)\le C
       &&\text{when }\eta=1.
\end{aligned}
\]
Letting $R\to\infty$ gives
\begin{equation}\label{eq:finite-dissipation}
 \int_0^\infty\int_{\Ct}\frac{G(X)}{hK}h_t^2\,\dd\sigma\,\dd t<\infty.
\end{equation}
Choose $t_j\in[j,j+1]$ such that
\[
 0\le\mathcal D(t_j)\le\int_j^{j+1}\mathcal D(t)\,\dd t\longrightarrow0.
\]
After passing to a subsequence, the uniform estimates imply
\[
 h(\cdot,t_j)\to h_\infty\quad\text{in every }C^k(\Ct),\qquad
 h_\infty\ge c,\quad A[h_\infty]\ge c\sigma,\quad
 (h_\infty)_\mu=\cot\theta\,h_\infty.
\]
For the corresponding choice of the time factor, set
\[
 \eta_\infty=
 \begin{cases}
 N[h_\infty]/D[h_\infty],&\eta=N/D,\\
 1,&\eta=1.
 \end{cases}
\]
Continuity of the integrals and~\eqref{eq:scalar} give
\[
\begin{aligned}
 \eta(t_j)&\to\eta_\infty>0,\\
 h_t(\cdot,t_j)&\to h_\infty
 \left(1-\frac{\eta_\infty fK_\infty}
                    {\phi(h_\infty)G(X_\infty)}\right)
 \quad\text{in every }C^k.
\end{aligned}
\]
Taking the limit of $\mathcal D(t_j)$ yields
\[
 0=\int_{\Ct}\frac{h_\infty G(X_\infty)}{K_\infty}
 \left(1-\frac{\eta_\infty fK_\infty}
                    {\phi(h_\infty)G(X_\infty)}\right)^2\dd\sigma.
\]
The weight is positive and the integrand is continuous, so
\[
 \eta_\infty^{-1}\phi(h_\infty)G(X_\infty)
                 \det A[h_\infty]=f,
\]
which is~\eqref{eq:problem} with $c=\eta_\infty^{-1}$, and $c=1$
when $\eta=1$. Evenness passes to the limit in the first case.

Let $\rho_j(\omega)=\rho(\omega,t_j)$ and let $\rho_\infty$ be the
radial function of the limit body. Since $X[h(t_j)]\to X[h_\infty]$
smoothly and $\langle X,\nu\rangle/|X|\ge c>0$, the radial maps
and their inverses converge, and hence $\rho_j\to\rho_\infty$ uniformly.
For $\eta=N/D$, the common bounds $c\le\rho_j,\rho_\infty\le C$
and the radial representation give
\[
\begin{aligned}
 &\left|V_{G,\pm}(\widehat\Sigma_{t_j})
           -V_{G,\pm}(\widehat\Sigma_\infty)\right|\\
 &\quad\le\int_{\Sp}
   \left|\int_{\rho_\infty(\omega)}^{\rho_j(\omega)}
                   G(s\omega)s^n\,\dd s\right|\dd\omega
 \le C\|\rho_j-\rho_\infty\|_{C^0}\longrightarrow0.
\end{aligned}
\]
Thus the preserved weighted volume has the same value at the limit.
For $\eta=1$, the bounds~\eqref{eq:barrier-bounds} hold for all times.
Finally,
\[
 X[h(t_j)]=\nabla h(t_j)+h(t_j)\nu
 \longrightarrow\nabla h_\infty+h_\infty\nu=X[h_\infty]
 \quad\text{in every }C^k,
\]
which proves the asserted hypersurface convergence.
\end{proof}

\begin{proof}[Proof of Corollary~\ref{cor:scales}]
Fix $s_0>0$ and write
$F_s(\xi)=s^n\phi(sl(\xi))G(s\xi)$.
By~\eqref{ass:limits}, there exist $0<s_-<s_0<s_+$ such that
\[
 \min_{\Ct}F_{s_-}>\max_{\Ct}f,
 \qquad \max_{\Ct}F_{s_+}<\min_{\Ct}f.
\]
Consequently,
\[
 s_-^n\phi(s_-l(\xi))G(s_-\xi)>f(\xi)
 >s_+^n\phi(s_+l(\xi))G(s_+\xi),\qquad\xi\in\Ct.
\]
Theorem~\ref{thm:barrier}, with initial cap $s_0\Ct$, gives the result.
\end{proof}

\subsection{Uniqueness}
Assume~\eqref{ass:U}. The scale inequality~\eqref{eq:dilation-inequality}
already implies~\eqref{ass:limits}, so existence with $c=1$ follows
from Corollary~\ref{cor:scales}. The same inequality gives a direct
comparison proof of uniqueness, extending the role of $p-q$ in
\cite[Section 5]{Gao}.

Let $h_1,h_2$ solve~\eqref{eq:problem} with $c=1$, and put
$q=h_1/h_2$, $a=\max_{\Ct}q$. The Robin condition implies
\[
 q_\mu=\frac{(h_1)_\mu h_2-h_1(h_2)_\mu}{h_2^2}=0.
\]
At an interior maximum, $\nabla q=0$ and $\nabla^2q\le0$.
At a boundary maximum, tangential differentiation of $q_\mu=0$ gives
\[
 q_\alpha=q_\mu=0,\qquad
 q_{\alpha\mu}=-\cot\theta\,q_\alpha=0,\qquad
 (q_{\alpha\beta})\le0,\quad q_{\mu\mu}\le0.
\]
Thus $\nabla q=0$ and $\nabla^2q\le0$ hold at any maximum point.
Since $h_1=qh_2$, at this point
\[
\begin{aligned}
 h_1&=ah_2,\qquad\nabla h_1=a\nabla h_2,
       \qquad X[h_1]=aX[h_2],\\
 A[h_1]&=aA[h_2]+h_2\nabla^2q\le aA[h_2],\\
 \det A[h_1]&\le a^n\det A[h_2].
\end{aligned}
\]
If $a>1$, using~\eqref{eq:dilation-inequality} at $\xi_0$ yields
\[
\begin{aligned}
 f&=\phi(h_1)G(X[h_1])\det A[h_1]\\
  &\le a^n\phi(ah_2)G(aX[h_2])\det A[h_2]\\
  &\le a^{-\delta_*}\phi(h_2)G(X[h_2])\det A[h_2]
   =a^{-\delta_*}f<f,
\end{aligned}
\]
a contradiction. Hence $a\le1$ and $h_1\le h_2$.
Interchanging $h_1,h_2$ gives $h_2\le h_1$, so $h_1=h_2$.

\subsection{Decay of the speed and convergence of the flow}
Let $h$ be the solution with $\eta=1$ and $h_0=s_0l$.
Condition~\eqref{ass:limits} supplies strict barriers on either side
of any given $s_0>0$. All estimates above therefore apply, and in
particular
\[
 r=\frac{fK}{\phi(h)G(X)}\ge r_0>0.
\]

\begin{proposition}\label{prop:decay}
Under~\eqref{ass:U}, there is $\omega>0$ such that, for $t\ge1$ and
integers $k\ge0$,
\begin{equation}\label{eq:rates}
 \begin{aligned}
 \|\partial_t\log h(t)\|_{C^0(\Ct)}
       &\le C e^{-\omega(t-1)},\\
 \|h(t)-h_*\|_{C^k(\Ct)}
       &\le C_k e^{-\omega(t-1)/(k+1)},\\
 \|X[h(t)]-X[h_*]\|_{C^k(\Ct)}
       &\le C_k' e^{-\omega(t-1)/(k+2)},
 \end{aligned}
\end{equation}
where $h_*$ is the unique stationary solution with $c=1$.
The constants are independent of $t$.
\end{proposition}
\begin{proof}
Set $w=\partial_t\log h=1-r$. By~\eqref{eq:log-linearization},
\[
\begin{aligned}
 w_t&=-r\partial_t\log r\\
 &=rhA^{ij}w_{ij}
   +r(2A^{ij}h_i+h\gamma_j)w_j+rE(h,X)w,\\
 w_\mu&=\partial_t\left(\frac{h_\mu}{h}\right)
       =\partial_t(\cot\theta)=0.
\end{aligned}
\]
Thus
\begin{equation}\label{eq:speed-equation}
 \begin{cases}
 w_t=rhA^{ij}w_{ij}
       +r(2A^{ij}h_i+h\gamma_j)w_j+rE(h,X)w
       &\text{in }\Ct\times(0,\infty),\\
 w_\mu=0&\text{on }\partial\Ct\times(0,\infty).
 \end{cases}
\end{equation}
The established estimates and~\eqref{ass:U} give
\[
 c\sigma^{ij}\le rhA^{ij}\le C\sigma^{ij},\qquad
 |r(2A^{ij}h_i+h\gamma_j)|\le C,\qquad
 rE(h,X)\le-\delta_*r_0=:-\omega<0.
\]
For $z=e^{\omega(t-1)}w$, we have
\[
 z_t=rhA^{ij}z_{ij}+r(2A^{ij}h_i+h\gamma_j)z_j
       +(rE(h,X)+\omega)z,\qquad z_\mu=0.
\]
Since $rE+\omega\le0$, the maximum principle gives
\begin{equation}\label{eq:speed-decay}
 \|w(t)\|_{C^0}\le\|w(1)\|_{C^0}e^{-\omega(t-1)},\qquad t\ge1.
\end{equation}
For $t_2>t_1\ge1$, use $h_t=hw$ to obtain
\[
\begin{aligned}
 \|h(t_2)-h(t_1)\|_{C^0}
 &\le\int_{t_1}^{t_2}\|h(s)w(s)\|_{C^0}\,\dd s\\
 &\le\frac{M\|w(1)\|_{C^0}}\omega e^{-\omega(t_1-1)}.
\end{aligned}
\]
Hence $h(t)$ converges uniformly. By Theorem~\ref{thm:barrier}
and the uniqueness proved above, its limit is the stationary
support function $h_*$. Letting
$t_2\to\infty$ gives $\|h(t)-h_*\|_{C^0}\le Ce^{-\omega(t-1)}$.
For $k\ge1$, interpolation and Proposition~\ref{prop:longtime} yield
\[
\begin{aligned}
 \|h(t)-h_*\|_{C^k}
 &\le C_k\|h(t)-h_*\|_{C^0}^{1/(k+1)}
       \|h(t)-h_*\|_{C^{k+1}}^{k/(k+1)}
       +C_k\|h(t)-h_*\|_{C^0}\\
 &\le C_k e^{-\omega(t-1)/(k+1)}.
\end{aligned}
\]
Finally,
\[
\begin{aligned}
 \|X[h(t)]-X[h_*]\|_{C^k}
 &=\|\nabla(h(t)-h_*)+(h(t)-h_*)\nu\|_{C^k}\\
 &\le C_k\|h(t)-h_*\|_{C^{k+1}}
 \le C_k'e^{-\omega(t-1)/(k+2)},
\end{aligned}
\]
proving~\eqref{eq:rates}.
\end{proof}

The uniqueness argument and Proposition~\ref{prop:decay} complete the
proof of Theorem~\ref{thm:unique}. The decay exponents above are
convenient interpolation rates and are not asserted to be optimal.
All geometric convergence statements are understood in inverse
capillary Gauss coordinates, or equivalently up to tangential
reparametrization.

\section{Power weights and nonhomogeneous examples}\label{sec:examples}

\subsection{The capillary \texorpdfstring{$L_p$}{Lp} dual Minkowski problem}
For $\phi(s)=s^{1-p}$ and $G(Y)=|Y|^{q-n-1}$,
\eqref{eq:problem} becomes
\begin{equation}\label{eq:power-equation}
 \det(\nabla^2h+h\sigma)
 =c^{-1}fh^{p-1}(h^2+|\nabla h|^2)^{(n+1-q)/2},
 \qquad h_\mu=\cot\theta\,h.
\end{equation}
The integral conditions reduce to $p,q>0$ for~\eqref{ass:plus}
and $p,q<0$ for~\eqref{ass:minus}, while
\[
 s^n\phi(sl(\xi))G(s\xi)
       =s^{q-p}l(\xi)^{1-p}|\xi|^{q-n-1},\qquad E=q-p.
\]
Thus Theorem~\ref{thm:even} provides a volume-preserving flow with a
stationary subsequence when $pq>0$ and the prescribed functions are even.
Theorem~\ref{thm:unique} gives the following dynamical conclusion
without a symmetry assumption.

\begin{corollary}\label{cor:power}
Let $n\ge2$, $\theta\in(0,\pi/2)$ and $p>q$. For every positive
$f\in C^\infty(\Ct)$ and every initial spherical cap $s_0\Ct$,
flow~\eqref{eq:flow} with
$\phi(s)=s^{1-p}$, $G(Y)=|Y|^{q-n-1}$ and $\eta=1$ exists for all
positive times and converges smoothly to the unique solution
of~\eqref{eq:power-equation} with $c=1$. The estimates
\eqref{eq:rates} hold with $\omega=(p-q)r_0$.
\end{corollary}
\begin{proof}
For $\phi(s)=s^{1-p}$ and $G(Y)=|Y|^{q-n-1}$,
\[
 \frac{s\phi'(s)}{\phi(s)}=1-p,\qquad
 \gradR\log G(Y)\cdot Y=q-n-1,
\]
and hence
\[
 E(s,Y)=n+(1-p)+(q-n-1)=q-p=-(p-q)<0.
\]
Thus~\eqref{ass:U} holds with $\delta_*=p-q$.
Theorem~\ref{thm:unique} gives existence and convergence, and
Proposition~\ref{prop:decay} gives~\eqref{eq:rates} with
$\omega=\delta_*r_0=(p-q)r_0$.
\end{proof}

The elliptic existence and uniqueness statement for $p>q$ is already
proved in~\cite{Gao,HHY}. Corollary~\ref{cor:power} concerns the
convergent evolution toward that solution. The elementary dilation
relation for power weights also explains the constant in
Theorem~\ref{thm:even}: replacing $h$ by $ah$ multiplies
$\phi(h)G(X)\det A$ by $a^{q-p}$. If $p\ne q$, a limit with constant
$c$ therefore becomes a solution with $c=1$ after replacing $h$ by
$c^{1/(q-p)}h$. This changes its weighted volume by the factor
$c^{q/(q-p)}$. When $p=q\ne0$, dilation leaves the equation invariant
and cannot remove $c$; no uniqueness of the volume-constrained shape
is asserted here. The cases $p=0$ or $q=0$ do not satisfy the integral
conditions, but are covered by Theorem~\ref{thm:unique} whenever $p>q$.

\subsection{Position weights beyond the power case}
The general results are not limited to a homogeneous prescription.
For example, fix $\tau\ne0$ and take
\begin{equation}\label{eq:nonhomogeneous}
 \phi(s)=s^{-1}\exp\!\left(\tfrac14\sin\log s\right),\qquad
 G(Y)=|Y|^{-n}\exp\!\left(
       \tau\frac{Y_1}{|Y|}+\tfrac14\sin\log|Y|\right).
\end{equation}
Since $Y_1/|Y|$ is homogeneous of degree zero,
\[
 \gradR\left(\frac{Y_1}{|Y|}\right)\cdot Y=0,\qquad
 \frac{s\phi'(s)}{\phi(s)}=-1+\tfrac14\cos\log s,\qquad
 \gradR\log G(Y)\cdot Y=-n+\tfrac14\cos\log|Y|.
\]
Therefore
\[
 E(s,Y)=-1+\tfrac14\cos\log s+\tfrac14\cos\log|Y|
       \le-\tfrac12.
\]
Theorem~\ref{thm:unique} applies to every positive smooth $f$.
Both weights are nonhomogeneous, and $G$ is nonradial and not even.

This example also distinguishes a general position weight from a
prescribed function of $(\xi,u)$. At the north point of the cap,
choose two small interior perturbations of a spherical support,
with the same value $h=s$ and opposite gradients $\pm\varepsilon E_1$.
They remain positive, strictly convex and Robin, and their positions
are $sE_{n+1}\pm\varepsilon E_1$. The ratio of the two values of $G$ is
\[
 \exp\!\left(\frac{2\tau\varepsilon}{\sqrt{s^2+\varepsilon^2}}\right)
       \ne1.
\]
Thus the same $(\xi,u)$ can give different coefficients, even after
multiplication by a fixed positive function of $\xi$. This dependence
is the reason for retaining the position derivatives in the curvature
argument rather than treating them as derivatives of a prescribed
function of $(\xi,u)$.

There are also nonhomogeneous even examples for
Theorem~\ref{thm:even}. For $pq>0$ and $\tau\ne0$, set
\[
 \phi(s)=s^{1-p}e^{\frac14\sin\log s},\qquad
 G(Y)=|Y|^{q-n-1}
       e^{\tau(Y_1/|Y|)^2+\frac14\sin\log|Y|}.
\]
The exponential factors are bounded above and below by positive
constants. Therefore~\eqref{ass:plus} holds when $p,q>0$, and
\eqref{ass:minus} when $p,q<0$, without either weight being homogeneous.

\begin{remark}
The assumptions above are sufficient conditions, not necessary ones.
For example, constant weights give the classical capillary Minkowski
equation, whose necessary horizontal moment conditions are
$\int_{\Ct}\xi_\alpha f\,\dd\sigma=0$ for $1\le\alpha\le n$;
see~\cite{MWWMinkowski}. Thus an existence statement for arbitrary
positive data cannot hold without some restriction on the weights
or the data. The hypotheses of the two flows address this issue in
different ways: evenness and a preserved volume in
Theorem~\ref{thm:even}, and strict scale barriers in
Theorem~\ref{thm:barrier}.
\end{remark}

\section*{Acknowledgements}

This paper is a natural continuation of the author's previous work
on the capillary $L_p$ dual Minkowski problem~\cite{Gao}.
Here we use curvature flows to study the more general dual
Orlicz--Minkowski problem for capillary hypersurfaces.
The main ideas, research questions, and overall approach of this
work were developed by the author. During the preparation of this
paper, AI tools (ChatGPT and DeepSeek) were used as technical
assistants to help refine the details of examples, locate relevant
references, and improve the clarity of the exposition.
All mathematical statements, proofs, and conclusions were
independently verified by the author.


\begin{thebibliography}{99}
\setlength{\itemsep}{0pt}\fontsize{9.5}{11.2}\selectfont

\bibitem{Ale38} A. D. Aleksandrov, \emph{On the theory of mixed volumes. III. Extensions of two theorems of Minkowski on convex polyhedra to arbitrary convex bodies}, Mat. Sb., {\bf 3} (1938) 27--46.

\bibitem{Ale39} A. D. Aleksandrov, \emph{On the surface area measure of convex bodies}, Mat. Sb., {\bf 6} (1939) 167--174.

\bibitem{BF.JDE.266-2019.7980} K. J. B\"{o}r\"{o}czky and F. Fodor, \emph{The {$L_p$} dual Minkowski problem for {$p>1$} and {$q>0$}}, J. Differential Equations, {\bf266} (2019) 7980--8033.

\bibitem{BHP.JDG.109-2018.411} K. J. B\"{o}r\"{o}czky, M. Henk and H. Pollehn, \emph{Subspace
  concentration of dual curvature measures of symmetric convex bodies}, J. Differential Geom., {\bf 109} (2018) 411--429.

\bibitem{BLYZ.JAMS.26-2013.831} K. J. B{\"o}r{\"o}czky, E. Lutwak, D. Yang and G. Zhang, \emph{The
  logarithmic Minkowski problem}, J. Amer. Math. Soc., {\bf 26} (2013) 831--852.

\bibitem{BIS2019} P. Bryan, M. N. Ivaki and J. Scheuer,
\emph{A unified flow approach to smooth, even {$L_p$}-Minkowski problems},
Anal. PDE, {\bf 12} (2019) 259--280.

\bibitem{BISOrlicz} P. Bryan, M. N. Ivaki and J. Scheuer,
\emph{Orlicz--Minkowski flows},
Calc. Var. Partial Differential Equations, {\bf 60} (2021) Art. 41.

\bibitem{CHZ.MA.373-2019.953} C. Chen, Y. Huang and Y. Zhao, \emph{Smooth solutions to the {$L_p$} dual Minkowski problem}, Math. Ann., {\bf 373} (2019) 953--976.

\bibitem{CCL} H. Chen, S. Chen and Q. R. Li,
\emph{Variations of a class of Monge--Amp\`ere-type functionals and their applications},
Anal. PDE, {\bf 14} (2021) 689--716.

\bibitem{CLFlow} H. Chen and Q. R. Li,
\emph{The {$L_p$} dual Minkowski problem and related parabolic flows},
J. Funct. Anal., {\bf 281} (2021) Art. 109139.

\bibitem{CLLX} L. Chen, Y. N. Liu, J. Lu and N. Xiang,
\emph{Existence of smooth even solutions to the dual Orlicz--Minkowski problem},
J. Geom. Anal., {\bf 32} (2022) Art. 40.

\bibitem{CTWX2022} L. Chen, Q. Tu, D. Wu and N. Xiang,
\emph{Anisotropic Gauss curvature flows and their associated dual Orlicz--Minkowski problems},
Proc. Roy. Soc. Edinburgh Sect. A, {\bf 152} (2022) 148--162.

\bibitem{CY76} S. Y. Cheng and S. T. Yau, \emph{On the regularity of the solution of the {$n$}-dimensional Minkowski problem}, Comm. Pure Appl. Math., {\bf 29} (1976) 495--516.

\bibitem{CW.Adv.205-2006.33} K. S. Chou and X. J. Wang, \emph{The {$L_p$}-Minkowski problem and the Minkowski problem in centroaffine geometry}, Adv. Math., {\bf 205} (2006) 33--83.

\bibitem{GaoEarlier} Y. Gao,
\emph{The {$L_p$} dual Minkowski problem for capillary hypersurfaces},
arXiv:2510.12804v2.

\bibitem{Gao} Y. Gao,
\emph{The capillary {$L_p$} dual Minkowski problem for $p>q$ in higher dimensions},
arXiv:2609.07158v2, 2026.

\bibitem{HP.Adv.323-2018.114} M. Henk and H. Pollehn, \emph{Necessary subspace concentration conditions for the even dual Minkowski problem}, Adv. Math., {\bf 323} (2018) 114--141.

\bibitem{HHI} J. R. Hu, Y. X. Hu and M. N. Ivaki,
\emph{Capillary {$L_p$} Minkowski flows},
arXiv:2509.06110, 2025.

\bibitem{HHY} J. R. Hu, Y. X. Hu and B. Yang,
\emph{Capillary John ellipsoid theorem with applications to capillary curvature problems},
arXiv:2603.27252v4, 2026.

\bibitem{HICurvature} Y. X. Hu and M. N. Ivaki,
\emph{Capillary curvature images},
arXiv:2505.12921.

\bibitem{HISChristoffel} Y. X. Hu, M. N. Ivaki and J. Scheuer,
\emph{Capillary Christoffel--Minkowski problem},
arXiv:2504.09320.

\bibitem{HLYZ.Acta.216-2016.325} Y. Huang, E. Lutwak, D. Yang and G. Zhang, \emph{Geometric measures in the dual Brunn-Minkowski theory and their associated Minkowski
  problems}, Acta Math., {\bf 216} (2016) 325--388.

\bibitem{HZ.Adv.332-2018.57} Y. Huang and Y. Zhao, \emph{On the {$L_p$} dual Minkowski problem}, Adv. Math., {\bf 332} (2018) 57--84.

\bibitem{HLYZ.DCG.33-2005.699} D. Hug, E. Lutwak, D. Yang and G. Zhang, \emph{On the {$L_p$} Minkowski problem for polytopes}, Discrete Comput. Geom., {\bf 33} (2005) 699--715.

\bibitem{LiLiu} G. H. Li and C. Y. Liu,
\emph{Capillary Orlicz--Minkowski flow in the upper half-space},
arXiv:2601.14659, 2026.

\bibitem{LSW.JEMSJ.22-2020.893} Q. R. Li, W. Sheng and X. J. Wang, \emph{Flow by Gauss curvature to the Aleksandrov and dual Minkowski problems}, J. Eur. Math. Soc. (JEMS), {\bf 22} (2020) 893--923.

\bibitem{LSYY2022} Q. R. Li, W. Sheng, D. Ye and C. Yi,
\emph{A flow approach to the Musielak--Orlicz--Gauss image problem},
Adv. Math., {\bf 403} (2022) Art. 108379.

\bibitem{Lieberman} G. M. Lieberman,
\emph{Second Order Parabolic Differential Equations},
World Scientific, River Edge, NJ, 1996; revised edition, 2005.

\bibitem{Lieberman2002} G. M. Lieberman,
\emph{Higher regularity for nonlinear oblique derivative problems in Lipschitz domains},
Ann. Scuola Norm. Sup. Pisa Cl. Sci. (5), {\bf 1} (2002) 111--151.

\bibitem{LiebermanErrata} G. M. Lieberman,
\emph{Errata for the revised edition of Second Order Parabolic Differential Equations}.

\bibitem{LiuLu} Y. N. Liu and J. Lu,
\emph{A flow method for the dual Orlicz--Minkowski problem},
Trans. Amer. Math. Soc., {\bf 373} (2020) 5833--5853.

\bibitem{Lut.JDG.38-1993.131} E. Lutwak, \emph{The Brunn-Minkowski-Firey theory. I. Mixed volumes and the Minkowski problem}, J. Differential Geom., {\bf 38} (1993) 131--150.

\bibitem{Lut.Adv.118-1996.244} E. Lutwak, \emph{The Brunn-Minkowski-Firey theory. II. Affine and geominimal surface areas}, Adv. Math., {\bf 118} (1996) 244--294.

\bibitem{LYZ.TAMS.356-2004.4359} E. Lutwak, D. Yang and G. Zhang, \emph{On the {$L_p$}-Minkowski
  problem}, Trans. Amer. Math. Soc., {\bf 356} (2004) 4359--4370.

\bibitem{LYZ} E. Lutwak, D. Yang and G. Zhang,
\emph{{$L_p$} dual curvature measures},
Adv. Math., {\bf 329} (2018) 85--132.

\bibitem{MWWMinkowski} X. Q. Mei, G. F. Wang and L. J. Weng,
\emph{The capillary Minkowski problem},
Adv. Math., {\bf 469} (2025) 110230.

\bibitem{MWWLp} X. Q. Mei, G. F. Wang and L. J. Weng,
\emph{The capillary {$L_p$}-Minkowski problem},
Trans. Amer. Math. Soc., to appear, doi:10.1090/tran/9866;
arXiv:2505.07746v3.

\bibitem{MWWFlow} X. Q. Mei, G. F. Wang and L. J. Weng,
\emph{The capillary Gauss curvature flow},
arXiv:2506.09840, 2025.

\bibitem{MWWChristoffel} X. Q. Mei, G. F. Wang and L. J. Weng,
\emph{The capillary Christoffel--Minkowski problem},
Calc. Var. Partial Differential Equations, {\bf 65} (2026) Art. 186.

\bibitem{MWWSolitons} X. Q. Mei, G. F. Wang and L. J. Weng,
\emph{Uniqueness of capillary Gauss solitons},
arXiv:2608.22494v2.

\bibitem{Nir53} L. Nirenberg, \emph{The Weyl and Minkowski problems in differential geometry in the large}, Comm. Pure Appl. Math., {\bf 6} (1953) 337--394.

\bibitem{WZOrlicz} X. D. Wang and B. C. Zhu,
\emph{The capillary Orlicz--Minkowski problem},
J. Funct. Anal., {\bf 291} (2026), no. 8, Art. 111574.

\bibitem{Zha.CVPDE.56-2017.18} Y.~Zhao, \emph{The dual Minkowski problem for negative indices}, Calc. Var. Partial Differential Equations, {\bf 56} (2017) Art. 18.

\bibitem{Zha.JDG.110-2018.543} Y. Zhao, \emph{Existence of solutions to the even dual Minkowski problem}, J. Differential Geom., {\bf 110} (2018) 543--572.

\end{thebibliography}
\end{document}